\documentclass{article}

\usepackage{amsmath}
\usepackage{amssymb}
\usepackage{latexsym}
\usepackage{amsxtra}
\usepackage{amscd}
\usepackage{theorem}
\usepackage{xypic}
\usepackage{graphicx}

\theoremstyle{plain}

{\theorembodyfont{\slshape}

        \newtheorem{thm}{Theorem}[section]
        \newtheorem{cor}[thm]{Corollary}
        \newtheorem{lem}[thm]{Lemma}
        \newtheorem{prop}[thm]{Proposition}

}

{\theorembodyfont{\rmfamily}

        \newtheorem{defn}[thm]{Definition}
        
        \newtheorem{rem}[thm]{Remark}

}

\renewcommand{\em}{\sl}

\newcommand{\proof}{{\bf Proof:\ }}
\newcommand{\Endproof}{\hspace*{\fill} $\Box$ \vspace{1ex} \noindent }

\makeatletter
\renewcommand{\subsection}{\@startsection{subsection}{2}%
        {\z@}{-3.25ex plus -1ex minus-.2ex}{-1em}{\bf}}
\makeatother

\newcommand{\ZZ}{\mathbb{Z}}
\newcommand{\QQ}{\mathbb{Q}}

\newcommand{\NN}{\mathbb{N}}
\newcommand{\RR}{\mathbb{R}}
\newcommand{\PP}{\mathbb{P}}
\renewcommand{\AA}{\mathbb{A}}

\newcommand{\U}{\mathcal{U}}

\newcommand{\OO}{\mathcal{O}}

\newcommand{\D}{\mathcal{D}}
\newcommand{\HH}{\mathcal{H}}

\newcommand{\Aut}{{\rm Aut}}

\newcommand{\Gal}{{\rm Gal}}

\newcommand{\spe}{{\rm sp}}

\newcommand{\m}{\mathfrak{m}}

\newcommand{\Spec}{{\rm Spec\,}}
\newcommand{\Spf}{{\rm Spf}}
\newcommand{\Spm}{{\rm Spm}}
\newcommand{\Frac}{{\rm Frac}}

\newcommand{\To}{\;\longrightarrow\;}
\newcommand{\iso}{\stackrel{\sim}{\to}}

\newcommand{\lpfeil}[1]{\stackrel{#1}{\To}}

\newcommand{\abs}[1]{\lvert#1\rvert}

\newcommand{\norm}[1]{\lVert#1\rVert}

\newcommand{\Xs}{\mathsf{X}}
\newcommand{\Us}{\mathsf{U}}

\newcommand{\As}{\mathsf{A}}
\newcommand{\Ys}{\mathsf{Y}}
\newcommand{\Zs}{\mathsf{Z}}
\newcommand{\Vs}{\mathsf{V}}

\newcommand{\Ds}{\mathsf{D}}

\newcommand{\Yt}{\tilde{Y}}
\newcommand{\Es}{\mathsf{E}}

\newcommand{\X}{\mathcal{X}}
\newcommand{\Y}{\mathcal{Y}}

\newcommand{\OOh}{\hat{\OO}}

\newcommand{\rig}{{^{\rm rig}}}
\newcommand{\an}{{^{\rm an}}}

\title{Another proof of the Semistable Reduction Theorem}

\author{Kai Arzdorf and Stefan Wewers}

\date{}

\begin{document}

\maketitle

\begin{abstract}
  We give a new proof of the Semistable Reduction Theorem for curves. The main
  idea is to present a curve $Y$ over a local field $K$ as a finite cover of
  the projective line $X=\PP^1_K$. By successive blowups (and after replacing
  $K$ by a suitable finite extension) we construct a semistable model of $X$
  whose normalization with respect to the cover is a semistable model of $Y$.
\end{abstract}

\section{Introduction}

\subsection{}

Let $K$ be a field which is complete with respect to a discrete valuation
$v$. We let $R$ denote the valuation ring of $v$, $\m\lhd R$
the maximal ideal of $R$ and $k:=R/\m$ the residue field. 

Let $X$ be a smooth projective and absolutely irreducible curve over $K$. A
{\em model} of $X$ is a normal, flat and proper $R$-scheme $X_R$ such that
$X_R\otimes_R K =X$. Given a model $X_R$ of $X$, its special fiber is denoted
by $X_s:=X_R\otimes_R k$. The $k$-scheme $X_s$ is proper, connected and of
pure dimension one. We say that $X_R$ is a {\em semistable model} of $X$ if
$X_s$ is {\em nodal}, i.e.\ all singular points are ordinary double points. We
say that $X$ has {\em semistable reduction} if it has a semistable model
$X_R$.

\begin{thm}[Semistable Reduction Theorem] \label{ssredthm} There exists a
  finite extension $L/K$ such that the curve $X_L:=X\otimes_K L$ has
  semistable reduction (w.r.t.\ the unique extension of $v$ to $L$).
\end{thm}

The first proof of this theorem was given by Deligne and Mumford
(\cite{DeligneMumford69}, Corollary 2.7). Since then, many more proofs have
appeared in the literature, see e.g.\ \cite{Abbes2000}.

\subsection{}

The question that originally motivated the present paper is: how can one
explicitly determine a semistable model of a given curve? In a way the proof
of Theorem \ref{ssredthm} by Deligne and Mumford is constructive: choose
$n\geq 3$ prime to the residue characteristic of $K$ and let $L/K$ be the
smallest field extension over which the $n$-torsion points of the jacobian of
$X$ become rational. Then the minimal regular model of $X_L$ is semistable. In
theory, this gives an algorithm to determine a semistable model. It seems,
however, that several steps in this algorithm are today still computationally
too expensive to be practical for curves of genus $g\geq 3$.

In this paper we work out a new proof of Theorem \ref{ssredthm} which we hope
will ultimately lead to a more practical algorithm. The starting point of our
investigation was a paper by M.\ Matignon (\cite{Matignon03}, see also
\cite{LehrMatignon06}), which gives an algorithm to compute the semistable
reduction of $p$-cyclic covers of the projective line (satisfying an
additional assumption). Trying to generalize Matignon's method to a more
general situation, we noticed that it could be used as a germ for a new proof
of Theorem \ref{ssredthm}. 

\subsection{}

Let us give a brief sketch of our proof. The first idea is to view the curve
under consideration as a finite cover of the projective line. So we start with
a smooth projective $K$-curve $Y$ and choose a nonconstant separable finite
morphism $\phi:Y\to X:=\PP^1_K$. For any model $X_R$ of $X$ we obtain a model
$Y_R$ of $Y$ and a finite $R$-morphism $\phi_R:Y_R\to X_R$ by normalization of
$X_R$ in $Y$. The goal is now to determine a semistable model $X_R$ of $X$
such that $Y_R$ is semistable as well. We show that this is possible after
replacing $K$ by a finite extension (Theorem \ref{relssredthm}) and
obtain Theorem \ref{ssredthm} as an immediate consequence.

In order to prove Theorem \ref{relssredthm}, we may assume that the cover
$\phi$ is Galois. Let $G$ denote the Galois group of the cover $\phi$. If the
order of $G$ is prime to the residue characteristic of $K$, then it is well
known how to obtain a model $X_R$ with the desired properties: it suffices to
take a semistable model which separates the branch points of the cover $\phi$
(see e.g.\ \cite{LiuAG}, \S 10.4). In particular, if the residue
characteristic is zero, then the Semistable Reduction Theorem is relatively
easy to prove.

In the general case, let $X_R$ be any semistable model of $X=\PP^1_K$. Let
$Y_R$ be the normalization of $X_R$ in $Y$. By a theorem of Epp (\cite{Epp73})
we may assume that the special fiber $Y_s$ of $Y_R$ is reduced. If $Y_R$ is
semistable then we are done. Otherwise, there exists a singular point $y\in
Y_s$ which is not an ordinary double point.  Let $x\in X_s$ denote the image
of $y$ under the map $\phi_R$. The crucial step of our proof is to show that
there exists a blowup $f:X_R'\to X_R$ with center $x$ which `improves the
situation'. To make this a bit more precise, let $Y_R'$ denote the
normalization of $X_R'$ in $Y$. Then the induced map $g:Y_R'\to Y_R$ is a
blowup with center $\phi_R^{-1}(x)$. We say that $X_R'$ is an {\em
  improvement} of $X_R$ at $y$ if the singularities of the special fiber of
$Y_R'$ which lie on the fiber $g^{-1}(y)$ are `less bad' than the singularity
$y\in Y_s$ (`badness' of singularities can be measured by a suitable numerical
invariant). Once the existence of an improvement has been shown, the proof of
Theorem \ref{relssredthm} is straightforward: start with some semistable model
$X_R$ of $X$ (e.g.\ the smooth model $\PP^1_R$) and repeatedly apply the above
improvement procedure. After a finite number of steps we obtain a model $X_R'$
whose normalization in $Y$ is a semistable model of $Y$.

Our proof that an improvement exists is local in the sense that it depends
only on the formal completions of $X_R$ at $x$ and of $Y_R$ at $y$. Instead of
working with formal schemes, the crucial step of the proof is phrased in the
language of rigid geometry, as follows.

Let $X\rig$ and $Y\rig$ denote the rigid analytic spaces associated to the
$K$-curves $X$ and $Y$. We consider the formal fiber $\Xs_x:=]x[_{X_R}\subset
X\rig$ of the point $x$, i.e.\ the subset of points on $X\rig$ which
specialize to $x$, and likewise the formal fiber $\Ys_y:=]y[_{Y_R}\subset
Y\rig$ of $y$. Then $\phi$ induces a finite Galois cover
$\phi_y:\Ys_y\to\Xs_x$ of smooth rigid analytic spaces of dimension one. The
fact that $x$ is a smooth point of the special fiber $X_s$ implies that
$\Xs_x$ is an open disk, i.e.\ isomorphic to the rigid space
\[
      \{\, t\in\AA_K^1 \,\mid\, \abs{t}<1 \,\}.
\]
Let $\Ds\subset\Xs_x$ be an affinoid disk. By this we mean that after a finite
extension of $K$ there exists a parameter $t$ as above and some
$\epsilon\in\sqrt{\abs{K^\times}}$, with $0<\epsilon<1$, such that $D$ is the
open subspace of $\Xs_x$ defined by the condition $\abs{t}\leq\epsilon$. To
$\Ds$ one can associate a blowup $f:X_R'\to X_R$ with center $x$ whose
exceptional divisor is a $(-1)$-curve. The blowup $g:Y_R'\to Y_R$ induced by
$f$ is associated to the affinoid subdomain
$\Us:=\phi_y^{-1}(\Ds)\subset\Ys_y$.

We say that the affinoid disk $\Ds\subset\Xs_x$ is {\em exhausting} if the
complement $\Ys_y-\phi_y^{-1}(\Ds)$ decomposes as a union of open annuli. Let
$\D$ denote the set of all exhausting disks. A
well known lemma (see e.g.\ \cite{BoschLuetkebohmert85}, Lemma 2.4) says that
every sufficiently large affinoid disk $\Ds\subset\Xs_s$ is exhausting. In
particular, $\D$ is nonempty. One easily shows:

\vspace{1ex}\noindent
{\bf Proposition:} {\em The blowup $f:X_R'\to X_R$ associated to $\Ds$ is an
  improvement at $y$ if and only if $\Ds$ is a minimal element of $\D$ (with
  respect to inclusion).}

\vspace{1ex} Therefore, we have reduced the proof of the Semistable Reduction
Theorem to the claim that the set $\D$ has a minimal element.\footnote{The
  second named author has learned the idea of reducing the semistable
  reduction theorem to the above statement from a lecture of Raynaud at a
  conference in Rennes in 2009 (\cite{RaynaudLecture1})}

Our proof of the existence of a minimal exhausting disk is divided into two
cases. We first assume that the Galois group $G$ of the cover is solvable. In
this case the proof can be easily reduced to the case that $G$ is cyclic of
prime order. Under the latter assumption, there is an explicit construction
of the minimal exhausting disk, based on methods introduced by Matignon in
\cite{Matignon03}. This is worked out in detail in the first author's thesis
\cite{KaiDiss}. 

If $G$ is not solvable then we argue by contradiction, and assume that the set
$\D$ of all exhausting disks does not have a minimum. To each disk $\Ds$ in
$\D$ we associate a point $x_{\Ds}\in\Xs\an$ on the Berkovich analytic space
associated to $\Xs$ (essentially, $x_{\Ds}$ corresponds to the maximum norm on
$\Ds$). A compactness argument shows that the sequence $x_{\Ds}$ converges to
a point $x_0\in\Xs\an$. Points on $\Xs\an$ fall into four different classes,
see \cite{Berkovich90}, \S 1.4. We then show that in each of these four cases
we can derive a contradiction, thus proving the claim.  A crucial fact used in
these arguments is that `inertia groups are solvable'.

\subsection{}

The argument sketched above is really quite different from the traditional
proofs and requires less heavy machinery than most of them. For
instance, we do not use \'etale cohomology nor the Picard functor nor
resolution of singularities. Our use of rigid analytic geometry is very
limited and could be easily replaced by more elementary arguments. Ultimately,
our proof relies on valuation theoretic arguments. In this sense it may be
considered to be similar in nature to Temkin's proof of the stable
modification theorem for families of curves (\cite{Temkin10}), although this
is a much deeper and more difficult result.


The solvable case of our proof is truly constructive and gives a concrete and
useful algorithm to compute semistable reduction of curves in the cases where
it applies. Examples where the curve is a cyclic cover of the projective line
of order $p$ are worked out in \cite{KaiDiss}. In the nonsolvable case our
argument is, as it is written down here, fundamentally
nonconstruction. Nevertheless we believe that a future variant will yield a
constructive and practical method as well.

\vspace{3ex}
\thanks{{\bf Acknowledgements:} We would like to thank Andrew Obus for useful
  comments on an earlier version of this paper.}

\section{Semistable reduction for covers}

\subsection{}

Let $K$ be a field which is complete with respect to a discrete valuation
$v:K^\times\to\RR$. We let $R$ denote the valuation ring of $v$, $\m\lhd R$
the maximal ideal of $R$ and $k:=R/\m$ the residue field. We also let $\pi$
denote a uniformizer of $R$; the particular choice of $\pi$ will play no
role. 

We make the additional assumption that the residue field $k$ is algebraically
closed. By \cite{LiuAG}, Lemma 10.4.5 this is no restriction of generality, as
far as the Semistable Reduction Theorem is concerned. 

We remark that our base ring $R$ is a complete discrete valuation ring and is
therefore excellent (see \cite{LiuAG}, \S 8.2). As a consequence, all schemes
and formal schemes occuring in this paper will be automatically
excellent. This fact will be used in several places throughout the paper. For
instance, if $A$ is a localization of an $R$-algebra of finite type, then $A$
is also excellent.

\subsection{Reduced special fiber and permanence} \label{permanence}

Let $X$ be a smooth projective and absolutely irreducible curve over $K$.  A
{\em pre-model} of $X$ is a flat and proper $R$-scheme $X_R$ such that
$X_R\otimes_R K=X$. Note that a pre-model is a model if and only if it is
normal.  

\begin{lem} \label{Serrelem}
  Let $X_R$ be a pre-model of $X$. Let $X_s:=X_R\otimes_R k$ be  the special
  fiber of $X_R$. Then $X_R$ is normal (i.e.\ a model) if and
  only if the following two conditions hold.
  \begin{itemize}
  \item[(a)]
    The special fiber $X_s$ has no embedded points (see \cite{LiuAG},
    Definition 7.1.6).
  \item[(b)] The local ring $\OO_{X_R,\eta}$ is a discrete valuation ring, for
    every generic point $\eta$ of $X_s$.
  \end{itemize}
\end{lem}

\proof This follows from Serre's criterion for normality, see \cite{LiuAG},
Theorem 2.23. Indeed, (b) is equivalent to the Condition (R1), whereas (a)
means that $X_s$ satisfies Condition (S1). By {\em loc.cit.}, Proposition
2.11, the latter is equivalent to Condition (S2) for $X_R$.
\Endproof

\begin{cor} \label{Serrecor}
  Let $X_R$ be a pre-model of $X$ with special fiber $X_s$.
  \begin{enumerate}
  \item
    If $X_s$ is reduced then $X_R$ is normal.
  \item
    Assume that $X_R$ is normal. Then $X_s$ is reduced if and only if $X_s$ is
    reduced in codimension zero (i.e.\  for every
    generic point $\eta\in X_s$ the local ring $\OO_{X_s,\eta}$ is a field).
  \end{enumerate}
\end{cor}

\proof Assume that $X_s$ is reduced. We have to show that Condition (a) and
(b) of the lemma hold true. This is obvious for Condition (a). In order to
verify Condition (b), let $\eta\in X_s$ be a generic point. Then
$A:=\OO_{X_R,\eta}$ is a noetherian local ring of dimension one. We have
\[
        \OO_{X_s,\eta}=A/\pi A.
\]
Therefore, our assumption that $X_s$ is reduced shows that $\pi A$ is the
maximal ideal of $A$. It follows that $A$ is a discrete valuation ring, i.e.\
Condition (b) of the lemma holds as well. This proves Assertion (i) of the
corollary. (See \cite{LiuAG}, Lemma 1.18, for a direct proof which does not
use the Serre criterion.)

For the proof of (ii) we assume that $X_R$ is normal. Then $X_s$ has no
embedded points (Condition (b)). Since $X_s$ has dimension one, this means
that $X_s$ is reduced if and only if it is reduced in codimension zero. 
\Endproof

Let $Y_R$ be a model of $Y$, and let $L/K$ be a finite
extension. Let $S$ denote the integral closure of $R$ in $L$. The {\em
  normalized base change} of $Y_R$ to $L$ is defined to be the normalization
$Y_S:=(Y_R\otimes_RS)^\sim$ of the scheme $Y_R\otimes_R S$. Clearly,
$Y_S$ is a model of $Y_L$.

\begin{prop} \label{Eppprop} Let $Y_R$ be a model of $Y$. Then there exists a
  finite extension $L/K$ such that the normalized base change
  $Y_S:=(Y_R\otimes_R S)^\sim$ of $Y_R$ to $L$ has a reduced special
  fiber. Furthermore, if $L'/L$ is any finite extension, then the usual base
  change $Y_S\otimes_S S'$ is normal. (Here $S$ and $S'$ denote the integral
  closures of $R$ in $L$ and $L'$.)
\end{prop}

\proof Let $\eta$ be a generic points of $Y_s$. Since $Y_R$ is assumed to be
normal, the local ring $A_\eta:=\OO_{Y_R,\eta}$ is a discrete valuation ring
dominating $R$, by Condition (b) of Lemma \ref{Serrelem}. By Corollary
\ref{Serrecor} (ii), the special fiber is reduced if and only if
\[
     \OO_{Y_s,\eta}=A_\eta/\pi A_\eta
\]
is a  field. The latter condition holds if and only if $\pi$ is a uniformizer
of $A_\eta$. 

Let $L/K$ be a finite extension, $S$ the valuation ring of $L$, $\pi'$ a
uniformizer of $S$ and $Y_S:=(Y_R\otimes_R S)^\sim$ the normalized base
change. Let $\eta'$ be a generic point of the special fiber of $Y_S$ lying
over $\eta$. Then the local ring $A_{\eta'}:=\OO_{Y_S,\eta'}$ is a discrete
valuation ring dominating $S$, . Moreover, $A_{\eta'}$ is a direct factor of the
integral closure of $A_\eta\otimes_R S$. Now it follows from a theorem of Epp
(\cite{Epp73}) that there exists a finite extension $L/K$ such that
$\pi'$ is a uniformizer for $A_\eta'$, for every generic point $\eta'$ of the
special fiber of $Y_S$. This implies that the special fiber of $Y_S$ is
reduced. 

Recall that the residue field $k$ of $K$ is assumed to be algebraically
closed. So every finite extension of $K$ has residue field $k$. In particular,
if $L'/L$ is a further finite extension, then $Y_S$ and $Y_S\otimes_S S'$ have
the same special fiber, which is reduced. Now it follows from Corollary
\ref{Serrecor} (i) that $Y_S\otimes_S S'$ is normal. This completes the proof
of the proposition. 
\Endproof  

\begin{defn} \label{permanencedef} 
  A model $X_R$ of $X$ with reduced special
  fiber is called {\em permanent}.\footnote{This terminology is also inspired
    by Raynaud's talk \cite{RaynaudLecture1}.}
\end{defn}

Proposition \ref{Eppprop} says that every given model of $X$ becomes permanent
after normalized base change to a suitable finite extension of the
base field. Furthermore, permanent models are permanent in the sense that
their special fibers are unchanged under any finite extension of the base
field. Therefore, we may always assume, while proving the Semistable Reduction
Theorem, that any given model is permanent.

\subsection{} \label{localinvariants}

Let $Y_R$ be a permanent model of $Y$ with special fiber
$Y_s$. Let $\Yt_s$ denote the normalization of $Y_s$. Note that $\Yt_s$ is a
smooth (not necessarily connected) $k$-curve and that we have a finite
morphism $p:\Yt_s\to Y_s$ which is an isomorphism when restricted to the
smooth part of $Y_s$. For a closed point $y\in Y_s$ we set
\[
       \delta_y:=\dim_k(p_*\OO_{\Yt_s}/\OO_{Y_s})_y
\]
and 
\[
       m_y:=\abs{p^{-1}(y)}.
\]
It is easy to see that $\delta_y\geq m_y-1$. 

\begin{prop} \label{localprop}
Let $y$ be closed point of $Y_s$.
\begin{enumerate}
\item The point $y$ is a smooth point of $Y_s$ if and only if
  $\delta_y=0$.
\item
  The point $y$ is an ordinary double point of $Y_s$ if and only if $\delta_y=1$
  and $m_y=2$. 
\item
  We have
  \begin{equation} \label{genuseq}
    g = 1+\sum_V \,(g_V-1) + \sum_{y\in Y_s} \delta_y.
  \end{equation}
  Here $V$ runs over the irreducible components of $\Yt_s$ and $g_V$ denotes the
  genus of the normalization of $V$.
\end{enumerate}
\end{prop}  

\proof
This is well known. See for e.g.\ \cite{LiuAG}, Proposition 7.5.4 and
Proposition 7.5.15.
\Endproof

\subsection{}

Let $Y_R$ be a permanent model of $Y$. Let $Y_s$ denote the special fiber of
$Y_R$. 

\begin{defn} \label{moddef}
  A {\em modification} of $Y_R$ is an $R$-morphism $f:Y_R'\to Y_R$, where
  $Y_R'$ is another model of $Y$  and
  $f$ is the identity on the general fiber. The modification $f$ is called
  {\em permanent} if $Y_R'$ is permanent. The subset of $Y_s$ where $f$
  is not an isomorphism is called the {\em center} of $f$. 
\end{defn}

Note that a modification $f:Y_R'\to Y_R$ has connected fibers because $Y_R$ is
normal. An irreducible component $W$ of the special fiber $Y_s'$ of $Y_R'$ is
called {\em exceptional} of $f(W)$ is a closed point of $Y_s$. The union of
the exceptional components is called the {\em exceptional divisor}. The union
of the irreducible components of $Y_s'$ which are not exceptional is called
the {\em strict transform} of $Y_s$. The normalization $p:\tilde{Y}_s\to Y_s$
factors through a finite map $p':\tilde{Y}_s\to Y_s'$. The image of $p'$ is
precisely the strict transform.

\begin{defn} \label{simpledef} 
  A permanent modification $f:Y_R'\to Y_R$ is
  called {\em simple} if the following holds:
  \begin{enumerate}
  \item
    Every exceptional component intersects the strict transform. 
  \item
    Every point of intersection of an exceptional component with the strict
    transform is an ordinary double point of $Y_s'$. 
  \end{enumerate}
\end{defn}

\begin{defn} \label{improvedef} Let $f:Y_R'\to Y_R$ be a simple modification
  and $y\in Y_s$ a singular point. Then $f$ is called an {\em
    improvement at $y$} if for every closed point $y'\in f^{-1}(y)$ which does
  not lie on the strict transform of $Y_s$ we have
  \[
          \delta_{y'}<\delta_y \qquad\text{or}\qquad 
             \delta_{y'}=\delta_y,\;m_{y'}>m_y.
  \]
\end{defn}

\begin{lem} \label{improvelem} Let $f:Y_R'\to Y_R$ be a simple modification
  and $y\in Y_s$ a singular point which lies in the center of
  $f$. Assume that $f$ is {\em not} an improvement at
  $y$. Then:
  \begin{enumerate}
  \item
    The fiber $W:=f^{-1}(y)$ has a unique singular point
    $y'$. 
  \item Every irreducible component of $W$ intersects the strict transform of
    $Y_s$ in a unique point distinct from $y'$.
  \item
   The normalization of every irreducible component of $W$ has genus zero and
    contains a unique point lying over $y'$.
  \end{enumerate}
\end{lem}

\proof Let $Z\subset Y_s'$ denote the strict transform of $Y_s$. Write
$W\cap Z=\{y_1',\ldots,y_m'\}$. By Condition (iii) of Definition
\ref{simpledef}, each point $y_i'$ is an ordinary double point of $Y_s'$ and
hence a smooth point of $Z$ and of $W$. It follows that the finite map
$p':\tilde{Y}_s\to Z$ induces a bijection between the fiber $p^{-1}(y)$ (where
$p:\tilde{Y}_s\to Y_s$ is the normalization) and the set $W\cap Z$. Hence
$m=m_y$, in the notation of \S \ref{localinvariants}. We also have
$\delta_{y_i'}=1$ by Proposition \ref{localprop} (ii).

Set $U:=W-Z$ and let $S$ denote the set of irreducible components of
$W$. Comparing the two expressions for the genus $g$ obtained by applying
Proposition \ref{localprop} (iii) to $Y_R$ and to $Y_R'$, one easily shows
that
\begin{equation} \label{improveeq1}
  \delta_y = m_y -\abs{S} + \sum_{V\in S} g_V + \sum_{y'\in U}
    \delta_{y'}.
\end{equation}
It follows from Condition (ii) of Definition \ref{simpledef} that every
component $V\in S$ contains at least one of the points $y_i'$ (and this
$y_i'$ lies on no other component in $S$). Therefore,
$\abs{S}\leq m_y$, and so \eqref{improveeq1} gives the inequality
\begin{equation} \label{improveeq2}
\delta_y \geq \sum_{V\in S} g_V + \sum_{y'\in U} \delta_{y'}.
\end{equation}  
Since by assumption $f$ is not an improvement at $y$, there exists at least
one point $y'\in U$ with $\delta_{y'}\geq\delta_y$. But then
\eqref{improveeq2} implies that $\delta_{y'}=\delta_y\geq 1$, $\delta_{y''}=0$
for all $y''\in U\backslash\{y'\}$ and $g_V=0$ for all $V\in S$. This prove
(i) and the first half of (iii). Our argument also shows that the inequality
\eqref{improveeq2} is actually an equality. It follows that $m_y=\abs{S}$, and
this proves (ii).  Since $W$ is connected and $y'$ the only singular point,
every irreducible component must pass through $y'$. This shows that
$\abs{S}\leq m_{y'}$.  Finally, our assumption that $f$ is not an improvement
implies $m_{y'}\leq m_y=\abs{S}$, which proves the second half of (iii).
\Endproof

\subsection{} \label{relative}

Let $K(Y)$ denote the function field of $Y$. The extension $K(Y)/K$ is a
regular extension of transcendence degree one. Therefore, there exists an
element $x\in K(Y)$ such that $K(Y)/K(x)$ is finite and separable. The choice
of $x$ corresponds to a finite separable morphism $\phi:Y\to X:=\PP^1_K$ (we
identify the rational function field $K(x)$ with the function field of
$\PP^1_K$).

We will prove the following `relative version' of the Semistable Reduction
Theorem. 

\begin{thm} \label{relssredthm}
  Let $\phi:Y\to X:=\PP^1_K$ be as above. Then (after replacing $K$ by a
  finite extension) there exists a semistable model $X_R$ of $X$
  such that the normalization $Y_R$ of $X_R$ in $K(Y)$ is a semistable model
  of $Y$. 
\end{thm}

Obviously, Theorem \ref{relssredthm} implies Theorem \ref{ssredthm}. The proof
of Theorem \ref{relssredthm} will occupy the rest of this paper. We start with
a preliminary remark.

\begin{prop} \label{relprop1} For the proof of Theorem \ref{relssredthm} we
  may assume that the extension $K(Y)/K(x)$ is Galois.
\end{prop} 

\proof Let $\tilde{Y}$ be the smooth projective curve whose function field
$K(\tilde{Y})$ is the Galois closure of $K(Y)/K(x)$. Let
$G:=\Gal(K(\tilde{Y})/K(x))$ denote the Galois group and $H\subset G$ the
subgroup corresponding to $K(Y)$. Then $G$ acts on $\tilde{Y}$ and the natural
map $\tilde{Y}\to Y$ identifies $Y$ with the quotient curve $\tilde{Y}/H$. 

Let $X_R$ be a semistable model of $X$ such that its normalization
$\tilde{Y}_R$ in $K(\tilde{Y})$ is a semistable model of $\tilde{Y}$. Then
$Y_R:=\tilde{Y}_R/H$ is a semistable model of $Y$, see \cite{LiuAG}, Proposition
10.3.48. Moreover, the map $\phi:Y\to X$ extends to a finite map $Y_R\to
X_R$. Since $Y_R$ is normal, $Y_R$ is the normalization of $X_R$ in
$K(Y)$. This proves the proposition.
\Endproof

\subsection{} \label{relative2}

We can now formulate our strategy to construct a semistable model of $Y$. We
choose a finite separable map $\phi:Y\to X:=\PP^1_K$. By Proposition
\ref{relprop1}, we may assume that $\phi$ is a Galois cover, with Galois group
$G$. Let $X_R$ be a semistable model of $X$, and let $Y_R$ denote the
normalization of $X_R$ in $Y$. Then $Y_R$ is a $G$-equivariant model of $Y$
such that $Y_R/G=X_R$. By Proposition \ref{Eppprop} we may assume that $X_R$
and $Y_R$ are permanent. Let $\phi_s:Y_s\to X_s$ be the finite map induced by
$\phi$. Note that $\phi_s$ is invariant under the induced action of $G$ on
$Y_s$ and that the induced map $Y_s/G\to X_s$ is a homeomorphism. However, it
may not be an isomorphism. It is an isomorphism only if $G$ acts faithfully on
each irreducible component of $Y_s$.

\begin{defn} \label{criticaldef} A closed point $x\in X_s$ is called {\em
    critical} with respect to $\phi$ if the inverse image $\phi_s^{-1}(x)$
  contains a non-nodal point of $Y_s$. We say that the semistable model $X_R$
  is {\em admissible} with respect to $\phi$ if every critical point $x$ is a
  smooth point of $X_s$.
\end{defn}

The following proposition is the crucial step in our proof of the Semistable
Reduction Theorem. 

\begin{prop} \label{relimproveprop} Let $X_R$ be an admissible semistable
  model of $X$, relative to $\phi$. Let $Y_R$ be the normalization of $X_R$ in
  $Y$ (which we assume is permanent). Let $x\in X_s$ be a critical point. Then
  (after replacing $K$ by a finite extension) there exists a simple
  modification $f:X_R'\to X_R$ with center $x$ such that the following holds.
  \begin{enumerate}
  \item Let $Y_R'$ denote the normalization of $X_R'$ in $Y$ (which we assume
    is permanent). Then the induced map $g:Y_R'\to Y_R$ is a simple
    modification.
  \item
    The modification $g$ is an improvement at every point
    $y\in\phi_s^{-1}(x)$. 
  \end{enumerate}
\end{prop}

The proof of this proposition is given in the remaining sections, starting
with \S \ref{rigid}.

\subsection{} \label{relative3}

Assuming Proposition \ref{relimproveprop} for the moment we can give a proof
of Theorem \ref{relssredthm}. Let $\phi:Y\to X=\PP^1_K$ be as above. Then the
smooth model $X_R:=\PP^1_R$ is clearly admissible with respect to $\phi$. Let
$Y_R$ be the normalization of $X_R$ in $Y$.

If $Y_R$ is a semistable model then we are done. Otherwise, there exists a
critical point $x\in X_s$. Since the inverse image $\phi_s^{-1}(x)$ is a
single $G$-orbit, the invariant $\delta_y$ defined in \S \ref{localinvariants}
is the same for all $y\in\phi_s^{-1}(x)$. Hence we may write
$\delta_x:=\delta_y$. 

Let $f:X_R'\to X_R$ be a simple modification as in Proposition
\ref{relimproveprop}, relative to $x$. Since $x$ is a smooth point of $X_s$,
the fiber $f^{-1}(x)\subset X_s'$ is a smooth curve of genus zero,
intersecting the strict transform of $X_s$ transversally in a unique point
$x_0'$ (this follows easily from \eqref{improveeq1}). In particular, the model
$X_R'$ is semistable.

If the normalization $Y_R'$ of $X_R'$ in $Y$ is semistable, then we are
done. Other\-wise, let $x'\in f^{-1}(x)$ be a critical point with respect to
$\phi$. It follows from Condition (ii) in Proposition \ref{relimproveprop} and
Definition \ref{simpledef} that $x'$ is a smooth point of $X_s$ (this shows
that the model $X_R'$ is admissible for $\phi$). By Definition
\ref{improvedef} we have $\delta_{x'}<\delta_x$ or $\delta_{x'}=\delta_x$,
$m_{x'}>m_x$. 

All in all we see that by repeated application of Proposition
\ref{relimproveprop} we can either strictly decrease the invariant $\delta_x$
or keep it constant and increase $m_x$. Since $\delta_x\geq 0$ and
$m_x\leq\delta_x+1$, this process has to stop after a finite number of
steps. It ends with a semistable model $X_R$ whose normalization in $Y$ is a
semistable model of $Y$.
\Endproof
 
\section{The rigid analytic point of view} \label{rigid}

\subsection{} \label{tubes}

We keep the assumption on our base field $K$. In the context of rigid analytic
geometry it is more convenient to work with an absolute value instead of with
an (exponential) valuation. We therefore choose a real constant $0<q<1$ and set
$\abs{a}:=q^{v(a)}$ for $a\in K$.  

Let $X$ be a smooth projective $K$-curve. We let $X\rig$ denote the rigid
analytic space associated to $X$, see e.g.\ \cite{BGR} or
\cite{FresnelvdPut}. Recall that the set underlying $X\rig$ is simply the
set of closed points of $X$. 

Let $X_R$ be a permanent model of $X$. Given a point $x\in X\rig$, its scheme
theoretic closure in $X_R$ intersects the special fiber $X_s$ in a unique
point $\bar{x}\in X_s$, called the {\em specialization} of $x$. The resulting
map
\[
      \spe_{X_R}:X\rig \to X_s
\]
is surjective and is called the {\em specialization map} of the model $X_R$.
 
Let $Z\subset X_s$ be a locally closed subscheme. Then the inverse image
\[
     ]Z[_{X_R}:=\spe_{X_R}^{-1}(Z)\subset X\rig
\]
is an open set in the $G$-topology for $X\rig$ and hence is a smooth rigid
analytic $K$-space. We call $]Z[_{X_R}$ the {\em tube} of $Z$ in $X_R$.
See e.g.\ \cite{Berthelot96}, \S 1. 

\begin{rem} \label{rigidrem} 
  Let $Z\subset X_s$ be a locally closed
  subscheme. Let $\X:=X_R|_Z\sphat$ be the formal completion of $X_R$ along
  $Z$.
\begin{enumerate}
\item   
  The tube $]Z[_{X_R}$ is canonically isomorphic to the {\em generic
    fiber} $\X_K$ of $\X$ as constructed in \cite{Berthelot96}, \S 1 (see also
  \cite{deJong95}, \S 7).
\item
  Let $\OO^\circ_{X\rig}$ denote the subsheaf of the structure sheaf on
  $X\rig$ consisting of functions that are bounded by $1$. Then we have a
  canonical isomorphism
  \[
          \Gamma(\X,\OO_{\X}) \iso \Gamma(]Z[_{X_R},\OO_{X\rig}^\circ),
  \]
  see \cite{deJong95}, Theorem 7.4.1. 
\item
  The definition of $]Z[_{X_R}$ is compatible with base change to any finite
  extension $L/K$. More precisely, we have a canonical isomorphism
  \[
        ]Z[_{X_R}\otimes_K L \cong ]Z[_{X_R\otimes_R S}.
  \]
  Here $S$ denotes the integral closure of $R$ in $L$, and we identify the
  special fiber of $X_R$ with that of $X_R\otimes_R S$. Note that
  $X_R\otimes_R S$ is again a permanent model by Proposition \ref{Eppprop}. 
\item The tube $]Z[_{X_R}$ is connected if and only if $Z$ is connected. This
  follows immediately from (ii): any idempotent function on
  $]Z[_{X_R}$ is analytic and bounded by $1$ and hence gives rise to an
  idempotent function on $Z=\X^{\rm red}$.
\item
  Combining (iii) and (iv) shows that the connected components of $]Z[_{X_R}$
  are absolutely connected.  
\item Suppose that $Z=\Spec(\bar{A})$ is an affine open subset of $X_s$. Then
  $\X=\Spf(A)$ is an affine formal scheme, where $A$ is flat and topologically
  of finite presentation over $R$, and is complete with respect to the
  $\pi$-adic topology (so $A$ is {\em admissible} in the terminology of
  \cite{BoschLuetkebohmert93}). Therefore, $A_K:=A\otimes_R K$ is an affinoid
  $K$-algebra.  Now it follows from the construction that
  $]Z[_{X_R}=\X_K=\Spm(A_K)$ is an affinoid subdomain of $X\rig$. 

  In this special case, (ii) says that 
  \[
     A=A_K^\circ:=\{\,f\in A_K \,\mid\, \norm{f}\leq 1 \,\}.
  \]
  Here $\norm{\,\cdot\,}$ denotes the maximum norm on the affinoid
  $\Spm(A_K)$. It follows that $\bar{A}=A/\pi A$ and that $Z=\Spec(\bar{A})$
  is the {\em canonical reduction} of the affinoid domain $]Z[_{X_R}$ (in the
  sense of \cite{FresnelvdPut}, \S 4.8).   
\end{enumerate}
\end{rem}

\subsection{} \label{residueclasses}

Let us now consider the case where $Z=\{x\}$ consists of a single closed point
of $X_s$. Then the tube $\Xs:=]x[_{X_R}$ is called the {\em residue class} of
$x$ (with respect to the model $X_R$). Let $A:=\OOh_{X_R,x}$ denote the
complete local ring of the model $X_R$ at $x$. By Remark \ref{rigidrem},
$\Xs$ can be identified with the generic fiber of the formal $R$-scheme
$\X=\Spf(A)$. Moreover, $A$ can be identified with the ring of analytic
functions on $\Xs$ bounded by $1$. It follows that the residue class
$\Xs$ depends, as a rigid analytic space, only on the completion of
$X_R$ at $x$.

\begin{defn} \label{boundarydef}
  An {\em open analytic curve} over $K$ is a rigid analytic $K$-space $\Xs$
  which becomes isomorphic, after a finite extension of $K$, to a residue
  class $]x[_{X_R}$, where $X_R$ is a permanent model of a smooth projective
  $K$-curve $X$ and $x\in X_s$ is a closed point of the special fiber.
  The formal $R$-scheme $\X=\Spf(A)$, where $A=\Gamma(\Xs,\OO_{\Xs}^\circ)$,
  is called the {\em canonical formal model} of $\Xs$.
  A {\em boundary point} of $\Xs$ is a generic point of $\Spec(A/\pi
  A)$. The set of boundary points of $\Xs$ is denoted by $\partial\Xs$.
\end{defn}

A boundary point $\eta\in\partial\Xs$ gives rise to a discrete valuation on
$\Frac(A)$. Its residue field $k(\eta)$ is a complete discrete valuation field
containing $k$. It is thus isomorphic to $k((t))$. 

Suppose that $\Xs=]x[_{X_R}$ is a residue class as above. Then
$A=\OOh_{X_R,x}$. It follows that a boundary point $\eta\in\partial\Xs$
corresponds to a local branch of $X_s$ through $x$. We obtain a natural
bijection between $\partial\Xs$ and the fiber $p^{-1}(x)$, where
$p:\tilde{X}_s\to X$ is the normalization of $X_s$. We write
$\bar{\eta}\in\tilde{X}_s$ for the point corresponding to
$\eta\in\partial\Xs$.

\begin{defn} \label{diskanndef}
  An open analytic curve $\Xs$ over $K$ is called an {\em open disk} if it is
  isomorphic to the standard open unit disk, i.e.\ to the rigid $K$-space
  \[
       \{\, t\in\AA^1_K \,\mid\, \abs{t}<1 \,\}.
  \]
  It is called an {\em open annulus} if it is isomorphic to 
  \[
        \{\, u\in\AA_K^1 \,\mid\, \epsilon<u<1 \,\},
  \]
  for some $\epsilon\in\sqrt{\abs{K^\times}}$, with $\epsilon<1$. 
\end{defn}

\begin{prop} \label{diskannprop}
  Let $X_R$ be a permanent model of a smooth projective curve over $K$. Let
  $x\in X_s$ be a closed point of the special fiber and let $\Xs:=]x[_{X_R}$
  denote the residue class of $x$. Then $x$ is a smooth point (resp.\ an
  ordinary double point) of $X_s$ if and only if $\Xs$ is an open
    disk (resp.\ an open annulus). 
\end{prop}

\proof (compare with \cite{BoschLuetkebohmert85}, Proposition 2.2 and 2.3)
Suppose first that $X_R=\PP^1_R$ is the projective line over $R$ and $x:=0\in
X_s=\PP^1_k$ the origin. It is then easy to see that the residue class of $x$
is the standard open unit disk, and that $\OOh_{X_R,x}=R[[t]]$, where $t$ is the
standard parameter on $\AA^1_R$.

Now let $X$ be any smooth $K$-curve, $X_R$ a permanent model and $x\in X_s$.
By the above, the residue class $]x[_{X_R}$ is an open disk if and only if
$\OOh_{X_R,x}\cong R[[t]]$. But the latter holds if and only if $x$ is a
smooth point. This proves the first equivalence.

The proof of the second equivalence is similar, but we have to be more careful
about the role of the base field $K$. As before, we can realize the standard
open annulus
\[
     \{\, u\in\AA_K^1 \,\mid\, \epsilon<u<1 \,\}
\]
as the residue class of a point $x$ on the special fiber of a model $X_R$ of
$\PP^1_K$. However, the model $X_R$ is permanent if and only if 
$\epsilon\in\abs{K^\times}$. If this is the case, then 
\[
    \OOh_{X_R,x} = R[[u,v\mid uv=a]],
\]
where $a\in K^\times$ is any element with $\abs{a}=\epsilon$. With this in mind,
the proof of the second equivalence is analogous to the proof of the first.
\Endproof 

\subsection{} \label{diskcover0}

We need a good notion of finite (Galois) covers of open analytic curves.

\begin{defn} \label{coverdef} Let $\Xs$ be an absolutely connected open
  analytic curve over $K$. An {\em admissible cover} of $\Xs$ is a finite and
  flat morphism $\phi:\Ys\to\Xs$ of rigid-analytic $K$-spaces, such that $\Ys$
  is also an open analytic curve. The cover $\phi$ is called a {\em regular
    Galois cover} if $\Ys$ is absolutely connected and the automorphism group
  $G:=\Aut(\phi)$ has order $\deg(\phi)$ (note that $\deg(\phi)\in\NN$ is well
  defined).
\end{defn}

Let $\phi:\Ys\to\Xs$ be an admissible cover. Let
$A:=\Gamma(\Xs,\OO_{\Xs}^\circ)$ and $B:=\Gamma(\Ys,\OO_{\Ys}^\circ)$. Then
$A$ is a normal and complete local domain with residue field $k$, and $B$ is a
finite $A$-algebra. Moreover, $B$ is a normal complete semilocal ring,
of the form $B=\oplus_i B_i$, where each $B_i$ is a normal complete local
domain with residue field $k$. After extending the base field $K$ we may
assume that $A/\pi A$ and $B/\pi B$ are reduced. Then $\phi$ is a regular
Galois cover if and only if $B$ is a domain and the field extension
$\Frac(B)/\Frac(A)$ is Galois. If this is the case, then the Galois group of
$\phi$ can be identified with the Galois group of $\Frac(B)/\Frac(A)$.

Let $\phi:\Ys\to\Xs$ be a regular Galois cover, with Galois group $G$. Let
$H\subset G$ be a subgroup. Then the quotient $\Zs:=\Ys/H$ is again an
absolutely irreducible open analytic curve. The induced maps $\Ys\to\Zs$ and
$\Zs\to\Xs$ are admissible covers. Moreover, $\Ys\to\Zs$ is a regular Galois
cover with Galois group $H$, and $\Zs\to\Xs$ is Galois if and only if $H$ is a
normal subgroup of $G$.

\subsection{}  \label{diskcover}

We now formulate our main result (Theorem \ref{diskthm}), and show that it
implies the Semistable Reduction Theorem. Let $\Xs$ be an open disk over $K$
and $\phi:\Ys\to\Xs$ a regular Galois cover. Let $G$ denote the Galois group
of $\phi$. (The assumption that $\Xs$ is a disk will be slightly relaxed in \S
\ref{solv}, but it will again be in force in \S \ref{nonsolv}).

\begin{defn} \label{exhaustdef}
\begin{enumerate}
\item
  A {\em parameter} for the open disk $\Xs$ is an element $t\in A$ such that
  $A=R[[t]]$.  
\item
  A subset $\Ds\subset\Xs$ is called a {\em closed disk} if there exists a
  parameter $t$ for the open disk $\Xs$ and a real number $\epsilon$,
  $0<\epsilon<1$, such that  
  \[
        \Ds = \Xs(\abs{t}\leq \epsilon).
  \]
  If $\Ds$ is also an affinoid subdomain then it is called an {\em affinoid
    disk} (this is the case iff $\epsilon\in\sqrt{\abs{K^\times}}$).   
\item
  An affinoid disk $\Ds\subset\Xs_x$ is called {\em exhausting} (with respect
  to $\phi$) if the complement $\Ys\backslash\phi^{-1}(\Ds)$ is the disjoint
  union of open annuli.  
\end{enumerate}
\end{defn}

We let $\D$ denote the set of all exhausting affinoid disks $\Ds\subset\Xs$.

\begin{lem} \label{BLlem} Let $t$ be a parameter for the open disk
  $\Xs$. Let $\epsilon,\epsilon'\in\sqrt{\abs{K^\times}}$ with
  $0<\epsilon<\epsilon'<1$.
  \begin{enumerate}
  \item
    If $X(\abs{t}\leq\epsilon)$ is exhausting then $X(\abs{t}\leq\epsilon')$
    is exhausting as well.
  \item There exists a constant $\epsilon_0<1$, such that
    $X(\abs{t}\leq\epsilon)$ is exhausting, for all $\epsilon\geq
    \epsilon_0$. 
  \end{enumerate} 
\end{lem}

\proof
This follows from \cite{BoschLuetkebohmert85}, Lemma 2.4.
\Endproof

\begin{cor}
  The set $\D$ is nonempty. Moreover, if $\Ds'\subset\Xs$ is an affinoid disk
  containing an element $\Ds\in\D$, then $\Ds'\in\D$. 
\end{cor}

The following theorem is a `local' version of Proposition \ref{relimproveprop}
and is really the main result of the present paper.

\begin{thm} \label{diskthm} Let $\phi:\Ys\to\Xs$ be a regular Galois cover of
  the open disk. Assume that $\Ys$ is {\em not} an open disk. Then the set
  $\D$ of all affinoid disks $\Ds\subset\Xs$ which are exhausting with respect
  to $\phi$ has a unique minimal element.
\end{thm}

The proof is given in \S \ref{solv} and \S \ref{nonsolv}, after some
preliminary remarks in \S \ref{prep}. In \S \ref{exhaust} we show that Theorem
\ref{diskthm} implies the Semistable Reduction Theorem.

\subsection{} \label{prep}

We fix a regular $G$-Galois cover $\phi:\Ys\to\Xs$ of the open disk. We let
$\X=\Spf(A)$ and $\Y=\Spf(B)$ denote the canonical formal models of $\Xs$ and
$\Ys$. We let $\eta\in\partial\Xs$ denote the unique boundary point of
$\Xs$. We let $k[\eta]\subset k(\eta)$ denote the valuation ring of the
residue field of $\eta$. In this section, we consider $\eta$ as a morphism of
formal schemes $\eta:\Spf(k[\eta])\to\X$. 

Let $\Ds\subset\Xs$ be an affinoid disk. It gives rise to a diagram of formal
$R$-schemes
\[\begin{CD}
      \Y' @>>> \Y \\
      @VVV     @VVV \\
      \X' @>>> \X, \\
\end{CD}\]
as follows. Let $t\in A$ be a parameter and $\epsilon\in\abs{K^\times}$ such
that $\Ds=\Xs(\abs{t}\leq\epsilon)$. Choose an element $a\in R$ with
$v(a)=\epsilon$ and let $\X'\to\X$ be the formal blowup of
the ideal $I:=(t,a)\lhd A$. Let $Z\subset\X'$ be the exceptional fiber (it is
equal to the reduced subscheme $(\X)^{\rm red}$, and it is isomorphic to
$\PP^1_k$). The morphism $\xi:\Spf(k[\eta])\to\X$ lifts uniquely to a
morphism $\xi':\Spf(k[\eta])\to\X'$.  Let $z\in Z$ denote the image of
$\xi'$ and $Z^\circ:=Z\backslash\{z\}$. Then
\[
      \Ds =]Z^\circ[_{\X'}.
\]

Let $\Y'$ be the normalization of the formal scheme $\X'$ in $\Ys$ (see
\cite{Conrad99}, \S 2.1). We call $\Y'$ the formal model of $\Ys$ induced by
$\Ds$. Let $W:=(\Y')^{\rm red}$ denote the reduced subscheme. Note that $W$ is
a connected projective $k$-curve. The canonical morphism $\Y'\to\X'$ restricts
to a finite map $W\to Z$. Let $\partial W\subset W$ denote the inverse image
of $z$ and $W^\circ:=W\backslash\partial W$. We have
\[
     \phi^{-1}(\Ds) = ]W^\circ[_{\Y'}.
\]
It follows that $\Ds$ is exhausting with respect to $\phi$ if and only if the
residue classes $]w[_{\Y'}$ are open annuli, for all $w\in\partial
W$. Actually, since $G$ acts transitively on the set $\partial W$, it suffices
that  this holds for one $w\in\partial W$.

Given a boundary point $\xi\in\partial\Ys$ the morphism
$\xi:\Spf(k[\xi])\to\Y$ lifts uniquely to $\xi':\Spf(k[\xi])\to\Y'$, and
the image $\bar{\xi}\in W$ of $\xi'$ lies in $\partial W$. We obtain a
surjective $G$-equivariant map
\begin{equation} \label{formalboundarymapeq}
    \partial\Ys\to\partial W.
\end{equation}
If $\Ds$ is exhausting, then this map is a bijection.
 
\begin{lem} \label{preplem}
  Assume that $\Ds\in\D$ is exhausting. Then the following holds.
  \begin{enumerate}
  \item
    The affinoid $\Es:=\phi^{-1}(\Ds)$ is absolutely connected.
  \item
    The cover $\Ys$ is an open disk if and only if $\Es$ is a closed disk.
  \item
    Assume that $\Ys$ is not an open disk. Then the set $\D$ of all exhausting
    closed disks is totally ordered with
    respect to inclusion. 
  \end{enumerate}
\end{lem}
  
\proof If $\Ds$ is exhausting, then $W^\circ\subset W$ is the complement of a
finite set of smooth points of $W$. Since $W$ is connected it follows that
$W^\circ$ is connected as well. By Remark \ref{rigidrem} (iv),(v) this implies
that $\phi^{-1}(\Ds)=]W^\circ[_{\Y'}$ is absolutely connected, proving (i).

The affinoid $\phi^{-1}(\Ds)=]W^\circ[_{\Y'}$ is a closed disk if and only if
$W^\circ\cong\AA^1_k$. Under the assumption that $\Ds$ is exhausting, this
holds if and only if $W\cong\PP^1_k$, $\Ys$ has a unique boundary point $\xi$
and the intersection of $W$ with the corresponding formal subscheme $\Spf
k[\xi]\subset\Y'$ is transversal. By Castelnuovo's criterion this holds if and
only if $\Y\cong\Spf R[[s]]$, i.e.\ $\Ys$ is an open disk. This proves (ii). 

To prove (iii) we let $\Ds'\subset\D$ be another exhausting disk, disjoint
from $\Ds$. We then have to show that $\Ys$ is an open disk. We may write
$\Ds=\{\,\abs{t}\leq\epsilon\}$, for some parameter $t$ and some
$\epsilon\in\abs{K^\times}$, $0<\epsilon<1$. Then there exists an
$\epsilon'>\epsilon$ sucht that $\D'$ is contained in a residue class of the
affinoid $\As:=\{\abs{t}=\epsilon'\}$ (a {\em closed annulus of thickness
  $0$}). Let $\Xs'\subset\As$ denote the residue class containing $\Ds'$. The
inverse image $\Ys':=\phi^{-1}(\Xs')$ is the disjoint union of residue classes
of the affinoid $\phi^{-1}(\As)$. Since $\Ds$ is exhausting, $\phi^{-1}(\As)$
is a disjoint union of closed annuli of thickness $0$ (this follows from
\cite{BoschLuetkebohmert85}, Lemma 2.4). In particular, $\phi^{-1}(\As)$ is an
affinoid with good reduction, and hence $\Ys'$ is the disjoint union of open
disks. Applying (ii) to the connected components of the cover $\Ys'\to\Xs'$ we
see that $\phi^{-1}(\Ds')$ is a disjoint union of closed disks. By (i),
$\phi^{-1}(\Ds')$ is connected and hence an open disk. Applying (ii) again
shows that $\Ys$ is an open disk, as desired. This finishes the proof of the
lemma.
\Endproof

\subsection{} \label{exhaust}

For the rest of this section we will show that Theorem \ref{diskthm} implies
Proposition \ref{relimproveprop} and hence the Semistable Reduction Theorem
(as explained in \S \ref{relative3}). 

We return to the situation considered in \S \ref{relative}. Let $Y$ be a
smooth projective and absolutely irreducible curve over $K$ and let $\phi:Y\to
X:=\PP^1_K$ be a finite separable and nonconstant morphism to the projective
line $X$. We assume that $\phi$ is a Galois cover, and let $G$ denote its
Galois group. The cover $\phi:Y\to X$ induces a finite morphism of rigid
analytic $K$-spaces $\phi\rig:Y\rig\to X\rig$. 

Let $X_R$ be a semistable model of $X$ which is admissible for $\phi$
(Definition \ref{criticaldef}). Let $Y_R$ denote the normalization of $X_R$ in
$Y$. We assume that $Y_R$ is a permanent model of $Y$. Let $x\in X_s$
be a critical point.  By assumption $x\in X_s$ is a smooth point. It follows
that the residue class $\Xs:=]x[_{X_R}$ is an open disk (Proposition
\ref{diskannprop}). Set $A:=\OOh_{X_R,x}$.

Let $\Ys:=\phi^{-1}(\Xs)\subset Y\rig$ denote the inverse image of the
residue class $\Xs$. This is an open rigid subspace of $Y\rig$ which is
invariant under the action of the Galois group $G$. In fact, 
\[
        \Ys = \bigcup_{y\in\phi_s^{-1}(x)} \Ys_y, \qquad \Ys_y:=]y[_{Y_R}
\]
is the disjoint union of the residue classes of the points $y\in Y_s$ lying
over $x$. The assumption that $x$ is a critical point is equivalent to the
statement that the residue classes $\Ys_y$ are not isomorphic to open disks.
Therefore, for any $y$ the induced cover $\phi_y:\Ys_y\to\Xs$ is a Galois
covers with Galois group $G_y={\rm Stab}_G(y)$ which satisfies the hypotheses
of Theorem \ref{diskthm}.

Let $\Ds\subset\Xs$ be an affinoid disk. After replacing $K$ by a finite
extension we may assume that $\Ds$ contains a $K$-rational point $P$. Choose
an element $t\in\OO_{X_R,x}$ which has a simple zero at $P$ and no other zero
on the residue class $\Xs$. Then $t$ is a parameter for the open disk
$\Xs$. Moreover, $\Ds=\Xs(\abs{t}\leq\epsilon)$ for some
$\epsilon\in\sqrt{\abs{K^\times}}$ with $0<\epsilon<1$. After a further
extension of $K$ we may assume that $\epsilon\in\abs{K^\times}$.

Choose an element $a\in R$ with $\abs{a}=\epsilon$, and let $f:X_R'\to X_R$ be
the blowup with center $x$ of the ideal $(t,a)\lhd\OO_{X_R,x}$. Let
$Z:=f^{-1}(x)$ be the exceptional divisor and $\X':=X_R'|_Z\sphat$ the formal
completion of $X_R'$ along $Z$. The natural morphism of formal $R$-schemes
$\X'\to\X=\Spf(A)$ is the formal blowup of the ideal $(t,a)\lhd A$, see
\cite{BoschLuetkebohmert93}, \S 2. By the explicit description of formal
blowups in {\em loc.cit.} one sees that $Z\cong\PP^1_k$ is a smooth curve of
genus zero which intersects the strict transform of $Y_s$ in a unique point
$z$ and that $z$ is an ordinary double point of $Y_s'$. In particular,
$f:X_R'\to X_R$ is a simple modification. Furthermore, 
\[
       \Ds = ]Z^\circ[_{X_R'}, \quad \text{with $Z^\circ:=Z-\{z\}$.}
\]

Let $g:Y_R'\to Y_R$ denote the modification induced by $f$ (i.e.\ $Y_R'$ is
the normalization of $X_R'$ in $Y$). After a finite extension of the
base field $K$ we may assume that the model $Y_R'$ is permanent. Let $W\subset
Y_s'$ denote the exceptional divisor of $g$. Since $Y_R'\to X_R'$ is finite,
it restricts to a finite map $W\to Z$. Moreover, the set $\partial W\subset W$
of points where $W$ intersects the strict transform of $Y_s'$ is precisely the
inverse image of $z$ in $W$. It follows that
\[
       \phi^{-1}(\Ds)=]W^\circ[_{Y_R'}, 
         \quad\text{with $W^\circ:=W\backslash\partial W$.}
\]
We say that $g$ is the modification of $Y_R$ induced by $\Ds$.  

Note that we have a natural surjective map
\begin{equation} \label{boundarymapeq}
   \partial\Ys \to \partial W
\end{equation}
mapping a boundary point $\xi\in\partial\Ys$ first to a point
$\bar{\xi}\in\tilde{Y}_s$ on the normalization of $Y_s$ (see \S
\ref{residueclasses}) and then to its image under the map $p':\tilde{Y}_s\to
Y_s'$ discussed after Definition \ref{moddef}. If the modification $g$ is
simple, then this map is a bijection.

We call the affinoid disk $\Ds$ {\em exhausting} if it is exhausting with
respect to the cover $\phi_y:\Ys_y\to\Xs$, in the sense of Definition
\ref{exhaustdef} and for some $y\in\phi^{-1}(x)$ (in fact this condition is
independent of $y$). As in \S \ref{diskcover} we let $\D$ denote the set of
all exhausting affinoid disks $\Ds\subset\Xs$.

\begin{figure}[ht]
  \centering
  \includegraphics[bb=0 0 900 500, scale=.22]{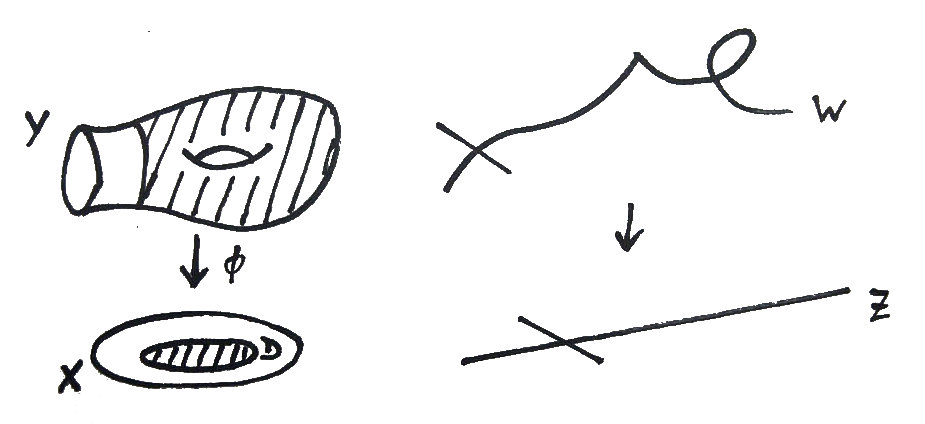}
  \caption{Here $\Ds$ is minimal exhausting.}
\end{figure}

Proposition \ref{relimproveprop} is now an immediate consequence of Theorem
\ref{diskthm} and the following lemma. 

\begin{lem} \label{disklem}
  Let $\Ds\subset\Xs$ be an affinoid disk and $g:Y_R'\to Y_R$ the induced
  modification.
  \begin{enumerate}
    \item
      The modification $g$ is simple if and only if $\Ds$ is exhausting. 
    \item Assume $\Ds$ is exhausting. Then $g:Y_R'\to Y_R$ is an improvement
      at every point $y\in\phi^{-1}(x)$ if and only if $\Ds$ is a minimal
      element of $\D$.
  \end{enumerate}
\end{lem}

\proof Every irreducible component of $W$ is a finite cover of $Z$ and
therefore contains a point $w\in\partial W$ lying above $z$. As remarked
above, these are precisely the points where $W$ intersects the strict
transform of $Y_s$. It follows that $g$ satisfies Condition (i) of Definition
\ref{simpledef}. Condition (ii) holds if and only if every point $w\in\partial
W$ is an ordinary double point of $Y_s'$. But the residue classes $]w[_{Y_R'}$
are precisely the connected components of
$\Ys\backslash\phi^{-1}(\Ds)$. Therefore, Condition (ii) of Definition
\ref{simpledef} holds if and only if $\Ds$ is exhausting (here we use
Proposition \ref{diskannprop}). This proves (i).

For the proof of (ii) we assume that $g$ is not an improvement at some
$y\in\phi^{-1}(x)$. Set $W_y:=g^{-1}(y)$; this is a connected component of
$W$. Set $W_y^\circ:=W_y\cap W^\circ$. By Lemma \ref{improvelem}, $W_y$ has a
unique singular point $y'$. All irreducible
components of $W_y$ have geometric genus zero and are smooth outside
$y'$. Moreover, they
intersect the strict transform in a unique point (which is an ordinary double
point distinct from $y'$) and have a unique branch passing through $y'$. 

Let $x'\in Z^\circ$ denote the image of $y'$. Then $x'$ is the only critical
point with respect to $\phi$ which lies on $Z$. Let $\Xs':=]x'[_{X_R'}$ denote
the residue class of $x'$ and let
$\Ys':=\phi^{-1}(\Xs)=]\phi^{-1}(x')[_{Y_R'}$ be the inverse image. Applying
Lemma \ref{BLlem} to the restriction $\phi':=\phi|_{\Ys'}:\Ys'\to\Xs'$ we
find an affinoid disk $\Ds'\subset\Xs'$ which is exhausting with respect to
$\phi'$. We claim that $\Ds'$, as an affinoid disk in $\Xs$, is 
exhausting with respect to $\phi$. Since $\Ds'$ is strictly contained in
$\Ds$, this claim would prove the `if' part of (ii).

To prove the claim, we consider the simple modification $g':Y_R''\to Y_R'$
induced from $\Ds$ (as an affinoid disk in $\Xs'$). The center of $g'$ is
precisely the singular locus of $W$. Let $W'\subset Y_s''$ denote the
strict transform  of $W$. By construction, 
\[
      \Xs\backslash\phi^{-1}(\Ds') = ]W'[_{Y_R''}.
\]
From the above description of $W$ and the fact that $g'$ is a simple
modification we see that $W'$ is isomorphic to the normalization of $W$. More
precisely, $W$ is a disjoint union of projective lines, each of which
intersects the rest of $Y_s''$ in exactly two points, which are ordinary
double points $Y_s''$. It follows that the tube $]W'[_{Y_R''}$ is the
disjoint union of open annuli. This proves the claim and hence the `if'-part
of (ii). The `only if' part is left to the reader (it is not used in the rest
of the paper).
\Endproof

\section{The solvable case} \label{solv}

In this section we prove the existence of a minimal exhausting disk (Theorem
\ref{diskthm}) under the assumption that the Galois group $G$ of the cover
$\phi:\Ys\to\Xs$ is solvable. The proof is by induction on the order of
$G$. The base case of the induction is when $G$ is cyclic of prime order, and
this case is treated in greater detail in \cite{KaiDiss}.  To make the
induction step work we actually have to consider a slightly more general
situation. Namely, we allow $\Xs$ to be either a disk or an annulus. It is
then natural to replace the notion of {\em exhausting affinoid disk} by {\em
  separating boundary domain}.

\subsection{} \label{solve0}

We keep all our assumptions on the base field $K$.  Let $\Xs$ be an open
analytic curve over $K$ which is either an open disk or an open annulus. Let
$A=\Gamma(\Xs,\OO_{\Xs}^\circ)$ be the ring of analytic functions on $\Xs$
bounded by $1$. Then $\Xs$ is the generic fiber of the formal $R$-scheme
$\X=\Spf(A)$, see \S \ref{tubes}-\S \ref{residueclasses}.


Let us choose a boundary point $\eta\in\partial\Xs$. If $\Xs$ is an annulus,
this amounts to choosing an `orientation' of $\Xs$ (if $\Xs$ is a disk, there is
no choice). A {\em parameter} for $\Xs$ (with respect to $\eta$) is an element
$t\in A$ which yields an isomorphism
\[
     \Xs \iso \{\, t\in\AA^1_K \mid \epsilon_0<\abs{t}<1 \,\},
\]
for some $\epsilon_0<0$ if $\Xs$ is a disk and with $\epsilon_0>0$
if $\Xs$ is an annulus. 
If $\Xs$ is
a disk, then this implies $A=R[[t]]$, and the new terminology agrees with the
old one. If $\Xs$ is an annulus, then there exists an element $a\in\m_R$ such
that $s:=a/t\in A$ is a parameter for $\Xs$ with respect to the boundary point
distinct from $\eta$, and we have $A=R[[t,s\mid ts=a]]$. 

By a {\em boundary domain} of $\Xs$ (containing $\eta$) we mean an open rigid
subspace $\Us\subset\Xs$ of the form
\[
       \Us = \Xs(\abs{t}>\epsilon),
\]
where $t$ is a parameter for $\Xs$ and $\epsilon\in\sqrt{\abs{K^\times}}$,
$\epsilon<1$.  So if $\Xs$ is an open disk, then $\Us=\Xs-\Ds$, where
$\Ds\subset\Xs$ is an affinoid disk. In any case, $\Us$ is an open annulus.

Let $\phi:\Ys\to\Xs$ be a regular Galois cover of $\Xs$, with Galois group
$G$. A boundary domain $\Us\subset\Xs$ is called {\em separating} (with
respect to $\phi$) if the inverse image $\phi^{-1}(\Us)$ is the disjoint union
of open annuli. Let $\U$ denote the set of all separating boundary domains. We
consider $\U$ as a partially ordered set by inclusion.

If $\Xs$ is an open disk, then a boundary domain $\Us\subset\Xs$ is separating
with respect to $\phi$ if and only if the affinoid disk
$\Ds:=\Xs\backslash\Us$ is exhausting with respect to $\phi$. It follows that
the set $\U$ has a unique maximum if and only if the set $\D$ has a
unique minimum. Therefore, the following proposition implies Theorem
\ref{diskthm} in case the Galois group $G$ is solvable. 

\begin{prop} \label{solvprop}
  Assume that
  \begin{itemize}
  \item[(a)] 
    The open analytic curve $\Ys$ is not an open disk. 
  \item[(b)]
    The Galois group $G$ of $\phi$ is solvable. 
  \end{itemize}
  Then the set $\U$ has a unique maximal element.
\end{prop}

For the proof of Proposition \ref{solvprop} we will use induction on the
order of $G$. The case where $G=\{1\}$ is trivial. Indeed, if $G=\{1\}$ then
$\Xs=\Ys$ is not a disk. It follows that $\Xs$ is an annulus, and that
$\Us:=\Xs$ is the unique maximal element of $\U$. After a preliminary argument
in \S \ref{solv1}, we prove the prime order case in \S \ref{solv2}--\S 
\ref{solv4}. Finally, the induction step is done in \S \ref{solv5}.

\subsection{} \label{solv1}

Let $x_1,\ldots,x_r\in\Xs$ be the pairwise distinct branch points of $\phi$
(which we assume to be $K$-rational). We claim that, in order to prove
Proposition \ref{solvprop}, we may assume that $r\leq 1$ if $\Xs$ is a disk
and $r=0$ if $\Xs$ is an annulus.

To prove this claim we assume that either $r\geq 2$ and $\Xs$ is a disk or
that $r\geq 1$ and $\Xs$ is an annulus. Under this condition there exists a
maximal boundary domain $\Us_0\subset\Xs$ containing none of the branch
points. Assume, moreover, that $\Us_0$ is not separating with respect to
$\phi$. Then $\phi^{-1}(\Us_0)$ is not a union of open disks. Furthermore, a
boundary domain $\Us\subset\Us_0$ is separating with respect to
$\phi:\Ys\to\Xs$ if and only it is separating with respect to the cover
$\phi^{-1}(\Us_0)\to\Us_0$. But $\phi^{-1}(\Us_0)\to\Us_0$ is \'etale by
choice of $\Us_0$. So in this case it suffices to prove the proposition under
the assumption that $\phi$ is \'etale (i.e.\ $r=0$).

Now consider the case that $\Us_0$ is separating. If it is maximal with this
property then we are done. Hence we may assume that $\Us_0$ is not
maximal. For simplicity, we also assume that $\Xs$ is a disk (the other case
is proved similarly). Then $\Ds_0:=\Xs\backslash\Us_0$ is a non-minimal
exhausting affinoid disk. In this situation it follows from Lemma
\ref{preplem} (iii) that there exists a unique residue class
$\Xs'\subset\Ds_0$ containing all exhausting disks strictly contained in
$\Ds_0$. Furthermore, for any closed affinoid disk $\Ds\subset\Xs'$, $\Ds$ is
exhausting with respect to $\phi$ if and only it is exhausting with respect to
the restricted cover $\phi^{-1}(\Xs')\to\Xs'$. But by the choice of $\Us_0$,
$\Ds_0$ is the smallest closed disk containing all $r\geq 2$ branch points. It
follows that the residue class $\Xs'\subset\Ds_0$ contains strictly less then
$r$ branch points. Our claim now follows by induction on the number $r$ of
branch points.

For the rest of the proof of Proposition \ref{solvprop} we may now assume that
$r\leq 1$, and that $\Xs$ is a disk if $r=1$.

\subsection{} \label{solv2}

Let us assume that the group $G$ is cyclic of prime order $\ell$. We have to
distinguish two main cases, of which the first is divided into two
subcases. We start with the assumption that the characteristic of $K$ is prime
to $\ell$. (For the time being, we make no assumption on the residue
characteristic of $K$.) Then, after replacing $K$ by some finite extension, we
may also assume that $K$ contains an $\ell$th rooth of unity. Moreover, the
cover $\phi:\Ys\to\Xs$ is given generically by a Kummer equation of the form
\[
      y^\ell = f,
\]
where $f\in A$ is not an $\ell$th power. More precisely, the ring
$B:=\OO(\Ys)^\circ$ contains an element $y$ which satisfies the above equation
and $B$ is the normalization of $A$ in the field $\Frac(A)[y]$. Clearly, the
ring $B$ is unchanged if we divide $f$ by an $\ell$th power in $A$. We may
therefore assume that the order of zero of $f$ at any point $x\in\Xs$ is
strictly less than $\ell$. Under this condition the zeroes of $f$ are
precisely the branch points of $\phi$.

We claim that there is no branch point, i.e.\ that $\phi$ is \'etale. To prove
this claim we assume the converse. Then $\Xs$ is a disk and there is exactly
one branch point, by our assumption made at the end of \S \ref{solv1}. We
choose a parameter $t$ for $\Xs$ such that the unique branch point is
$t=0$. Now the Weierstrass preparation theorem shows that $f$ is of the form
$f=c(t+a_2t^2+\ldots)\in A=R[[t]]$, with $c\neq 0$.  After a finite extension
of $K$, $c$ is an $\ell$th power, and we may therefore assume that
$c=1$. Applying \cite{KaiDiss}, Lemma 1.28, one shows that $B=A[y]$ and that
$\Ys$ is an open disk with parameter $y$. But this contradicts Assumption (a)
of Proposition \ref{solvprop}, proving the claim.

Let us now make the additional assumption that $\ell\neq p$, i.e.\ that the
residue characteristic of $K$ is prime to $\ell$. Let us choose a parameter
$t$ for $\Xs$. If $\Xs$ were a disk then $A=R[[t]]$ and, by the above claim,
$f=1+a_1t+a_2t^2+\ldots$. But then Hensel's
lemma shows that $f$ is an $\ell$th power in $A$. This contradicts our
assumption that $\phi$ is a regular Galois cover. We conclude that $\Xs$ is an
open annulus and hence $A=R[[t,s\mid ts=a]]]$. Using again Hensel's Lemma we
see that $f=ct^mu$, with $c\in R$, $m\geq 1$ and $u$ a unit of $A$ with
constant coefficient $1$. As before we may assume that $c=1$. Dividing $f$ by
a suitable power of $t^\ell$ we may also assume that $m<\ell$. In this
situation \cite{KaiDiss}, Lemma 1.31, shows that $\Ys$ is an open
annulus. This means that $\Xs$ itself is the maximal separating boundary
domain we are looking for. Therefore, Proposition \ref{solvprop} is proved in
the case $\abs{G}=\ell\neq p$.

\subsection{} \label{solv3}

We continue with the notation introduced in \S \ref{solv2}, but we now assume
that $\ell=p$ is the (positive) characteristic of the residue field $k$ of
$K$. (We keep the assumption that $\ell=p$ is prime to the characteristic of
$K$, but this now amounts to saying that $K$ has characteristic zero.) In this
case, Proposition \ref{solvprop} is proved in \cite{KaiDiss} (Theorem 2.1 and
Theorem 4.2). We briefly sketch the proof. For simplicity, we restrict to the
case where $\Xs$ is an open disk. The proof in the case where $\Xs$ is an
annulus uses the same methods, but is slightly more complicated.  

As we have seen before, we may assume that the cover $\phi$ is given
generically by an equation of the form
\[
     y^p= f, 
\]
where $f$ is a principal unit in $A=R[[t]]$. Note that Hensel's lemma is not
applicable anymore, and we cannot conclude that $f$ is a $p$th power. 

The first step is to compute the ring $B$. Let
$v_0:\Frac(A)\to\QQ\cup\{\infty\}$ denote the discrete valuation corresponding
to the prime ideal $\pi A\lhd A$ (normalized such that $v_0(p)=1$ which
implies that $v_0|K=v$). Since $\bar{A}:=A/\pi A$ is isomorphic to a power
series ring $k[[t]]$, the residue field of $v_0$ carries a canonical discrete
valuation $\bar{v}_0$ (normalized such that $\bar{v}_0(t)=1$). Let
\[
     v_\eta:A\to(\QQ\times\ZZ)\cup\{\infty\})
\]
denote the discrete valuation of rank two obtained as the composition of $v_0$
with $\bar{v}_0$. Here we consider the target set as an ordered group with
respect to the lexicographic ordering. More explicitly: if $t$ is a parameter
for $\Xs$ and $g=\sum_{i=0}^\infty g_it^i\in A$ then $v_\eta(g)=(\mu,m)$,
where
\[
     \mu=\min_iv(g_i), \qquad m=\min\{\, i \mid v(g_i)=\mu\}.
\]

\begin{prop} \label{pcyclicprop1}
  The maximum 
  \[
        (\mu,m) :=\max\{v_\eta(f-h^p) \mid h\in A^\times\}
  \]
  exists. Furthermore:
  \begin{enumerate}
  \item
    $0\leq\mu<p/(p-1)$ and $\mu/p\in v(K^\times)$.
  \item
    $m>1$ and $(m,p)=1$.
  \item  
    Choose $c\in K^\times$ such that $v(c)=\mu/p$ and $h\in A^\times$ such that
    $v_\eta(f-h^p)=(\mu,m)$. Then $B=A[w]$, where $w:=(y-h)/c$.
  \end{enumerate}
\end{prop}

\proof Let $h\in A^\times$ be given and set $(\mu,m):=v_\eta(f-h^p)$. Then
$\mu=v_0(f-h^p)\leq p/(p-1)$ (otherwise, completness of $A$ with respect to
$\pi A$ would show that $f\in A^p$). Since $\mu\in v(K^\times)$ takes values
in a discrete group, we may assume that $\mu$ takes the maximal possible
value. Let $\tilde{v}_0$ be an extension of $v_0$ to $\Frac(B)$; it
corresponds to a minimal prime ideal of $\bar{B}:=B/\pi B$. Our running
assumption says that $\bar{B}$ is reduced, which means that $\pi$ is a prime
element for $\tilde{v}_0$. A simple calculation using $\mu\leq p/(p-1)$ and
the equation $y^p=f$ shows that $p\tilde{v}_0(y-h)=v_0(f-h^p)=\mu$. It follows
that $\mu/p\in v(K^\times)$. Choose an element $c\in K$ with $v(c)=\mu/p$ and
set $g_0:=(f-h^p)/c^p\in A$. Then $w:=(y-h)/c$ satiesfies an irreducible
equation over $A$:
\begin{equation} \label{solvprop1eq1}
   w^p +\ldots+pc^{1-p}h^{p-1}w +g_0 = 0.
\end{equation}
It follows that $w\in B$. If $\mu=p/(p-1)$ then Hensel's lemma, applied to the
complete local ring $A$, would show that this equation is reducible. We
conclude that $\mu<p/(p-1)$. Now Part (i) of the proposition is proved. 

Let $\bar{w}$ denote the image of $w$ in $\bar{B}:=B/\pi B$ and $\bar{g}_0$
the image of $g_0$ in $\bar{A}:=A/\pi A=k[[t]]$. Then
$m=\bar{v}_0(\bar{g}_0)$. It is easy to see that for any choice of $h$ such
that $v_0(f-h^p)=\mu$ we have $\bar{g}_0\not\in\bar{A}^p$ (otherwise $\mu$
wouldn't be maximal). Moreover, choosing a different $h$ results in adding to
$\bar{g}_0$ an element of $\bar{A}^p$. We may therefore assume that
$(p,m)=1$. Furthermore, $m$ is now the maximal possible value for
$\bar{v}_0(\bar{g}_0)$. It follows that $(\mu,m)$ is the maximal possible
value for $v_\eta(f-h^p)$. This proves the first claim of the proposition and
Part (ii), except for the statement that $m>1$.

From \eqref{solvprop1eq1} we see that $\bar{w}$ satisfies the 
equation
\[
     \bar{w}^p = \bar{g}_0,
\]
which is irreducible because $\bar{g}_0\not\in\bar{A}^p$. We conclude that
$\tilde{v}_0$ is the unique extension of $v_0$ to $\Frac(B)$, with residue
field extension purely inseparabel of degree $p$. Using Serre's Normality
Criterion (as in Lemma \ref{Serrelem}) we also see that $A[w]$ is normal and
hence $B=A[w]$. This proves (iii). Moreover, if $m=1$ then we would have
$B=R[[w]]$, i.e.\ $\Ys$ was an open disk. This contradicts our assumption. Now
the proof of the proposition is complete.
\Endproof

We continue to use the notation $(\mu,m)$ from the above proposition. We note
that the values of $\mu$ and $m$ and do not change if we replace $K$ by any
finite extension. An element $h\in A$ such that $v_\eta(f-h^p)=(\mu,m)$ is
called a {\em best approximation} of $f$ (with respect to $v_\eta$) (cf.\
\cite{KaiDiss}, \S 2.2.1).  

\begin{lem}[$p$-Taylor expansion] \label{pcycliclem} Let $n\geq 1$ be given,
  and set
  \[
        \nu_n:=1+1/p+\ldots+1/p^n.
  \]
  After replacing $K$ by a finite extension, there
  exists a parameter $t$ for the disk $\Xs$ and an element $h\in A=R[[t]]$ 
  such that
  \[
         f-h^p=\sum_{j=1}^\infty a_j t^j 
  \]
  with $a_j\in\m_R$ and 
  \begin{equation} \label{pteq1}
       v(a_{pj})\geq \nu_n,
  \end{equation}
  for all $j$. 
\end{lem}
Following \cite{Matignon03}, we call $(t,h,a_j)$ a {\em $p$-Taylor expansion}
of $f$ of level $n$. 

\proof
See \cite{KaiDiss}, Proposition 2.12.
\Endproof

Since $\nu_n\to p/(p-1)$ for $n\to\infty$, we may choose $n$ such that
\begin{equation} \label{pteq2}
    \nu_n>\frac{p}{p-1}-\frac{p/(p-1)-\mu}{m}.
\end{equation}
Let $(t,h,a_j)$ be a $p$-Taylor expansion of $f$ of level $n$. Since
$\nu_n>\mu\geq v_0(f-h^p)$, it follows from \eqref{pteq2} that the minimum of
the valuations $v(a_j)$ occurs for an index $j$ which is prime to $p$. By
inspection of the proof of Proposition \ref{pcyclicprop1} one concludes that
$v_\eta(f-h^p)=(\mu,m)$, i.e.\ a $p$-Taylor expansion of $f$ yields a best
approximation (see \cite{KaiDiss}, Corollary 2.16).

We define 
\[
   \rho:=\min\,\{\,\frac{p/(p-1)-\mu}{m},\,\frac{v(a_j)-\mu}{m-j} \mid
                         1\leq j<m\,\}.
\]
Note that $\mu+m\rho\leq p/(p-1)$. If equality holds, we set
$k:=0$. Otherwise, we let $k$ denote the smallest index such that $1\leq k<m$
and $\rho=(v(a_k)-\mu)/(m-k)$. Then by definition the Newton polygon of
$f-h^p$ has a line segment of slope $-\rho$ over the intervall
$[k,\ldots,m]$. Moreover, it follows from \eqref{pteq1} and \eqref{pteq2} that
$(k,p)=1$.

One can check that $\rho$ and $k$ do not depend on the choice
of $h$ (see \cite{KaiDiss}, Proposition 2.28). However, $\rho$ and $k$ may
depend on the choice of the parameter $t$!

\begin{lem}
  After a finite extension of $K$, there exists a $p$-Taylor expansion
  $(t,h,a_j)$ of level $n$ such that $k\neq 1$.
\end{lem}

\proof
See \cite{KaiDiss}, Proposition 2.31. The idea is to use a `generic'
$p$-Taylor expansion 
\[
    f(t+T)-H^p=\sum_{j=1}^\infty A_jt^j,
\]
where the $A_j$ and the $t$-coefficients of $H$ are algebraic functions in
$T$. One has to show that, after a finite extension of $K$, there is a point
$T=\xi$ such that $A_1(\xi)=0$. Replacing the parameter $t$ by $t':=t-\xi$
then gives a $p$-Taylor expansion $(t',h',a_j')$ with $a_1'=0$. With respect
to this $p$-Taylor expansion we have $k\neq 1$. 
\Endproof

Now the following proposition completes the proof of Proposition
\ref{solvprop} in the special case considered in this subsection.

\begin{prop} \label{pcyclicprop2} Let $(t,h,a_j)$ be a $p$-Taylor expansion of
  $f$ such that $k\neq 1$ (notation as above).  Then
  $\Ds:=\Xs(\abs{t}\leq\rho)$ is the minimal exhausting disk with respect to
  $\phi$.
\end{prop}
    
\proof Let $\Us:=\Xs\backslash\Ds$. Choose an element $b\in R$ such that
$v(b)=\rho$ and set $t_1:=t/b$, $s_1:=t_1^{-1}$ (in general, this requires
again a finite extension of $K$). Then
\[
  A':=\Gamma(\Ds,\OO_{\Xs}^\circ) = R\{\{t_1\}\}
\]
(the ring of convergent power series over $R$) and
$B':=\Gamma(\phi^{-1}(\Ds),\OO_{\Ys}^\circ)$ is the integral closure of $A'$
in the extension of fraction fields given by $y^p=f$. Similarly,
\[
  A'':=\Gamma(\Us,\OO_{\Xs}^\circ) = R[[t,s_1\mid ts_1=b]],
\]
and $B'':=\Gamma(\phi^{-1}(\Us),\OO_{\Ys}^\circ)$ is the integral closure of
$A''$.

By construction, the Newton polygon of the power series $f-h^p=\sum_{j\geq
  1}a_jt^j$ has a line segment of slope $-\rho$ over the interval $[k,m]$ and,
moreover, $-\rho$ is the largest negative slope that occurs. It follows that
we can write $f-h^p$, as an element of the ring $A''$, in the form
\[
     f-h^p= c_1^pt^mu_1,
\]
where $c_1\in R$ is an element of valuation $\mu/p$ and $u_1$ is a principal
unit.  It follows that the element $w_1:=(y-h)/c_1$ is an element of $B''$
satisfying an irreducible equation over $A''$ of the form
\[
      w_1^p+\ldots+pc_1^{1-p}h^{p-1}=t^mu_1.
\]
Now it follows from \cite{KaiDiss}, Lemma 1.31, that $\phi^{-1}(\Us)$ is an
open annulus. We have shown that $\Ds$ is exhausting.

To show that $\Ds$ is the minimal disk with this property we shall provide an
explicit description of the canonical reduction $W^\circ$ of $\phi^{-1}(\Ds)$
(we freely use the notation set up in \S \ref{prep}) . Set
$\lambda:=\mu+m\rho$. Then $\lambda\leq p/(p-1)$ and we can write
\[
     f -h^p=c_2^pg,
\]
with $g\in A''$ and where $c_2\in R$ is chosen such that
$v(c_2)=\lambda/p$. Then the element $w:=(y-h)/c_2$ satisfies the integral
equation
\[
     w^p+\ldots+pc_2^{1-p}h^{p-1}= g
\]
over $A''$ and therefore $w\in B''$. Let $\bar{g}\in \bar{A}'':=A''/\pi
A''=k[t_1]$ (resp.\ $\bar{w}\in\bar{B}'':=B''/\pi B''$) denote the image
of $g$ (resp.\ of $w$). 

Suppose first that $\lambda=p/(p-1)$. Then $\bar{w}$ satisfies an
Artin-Schreier equation
\[
   \bar{w}^p+\bar{c}\bar{w}=\bar{g}.
\] 
Moreover, $\bar{g}\in k[t_1]$ is a polynomial in $t_1$ of degree $m>1$,
$(m,p)=1$. It follows that $W^\circ=\Spec k[t_1,\bar{w}]$ is a smooth and
connected affine curve of genus $g=(p-1)(m-1)/2>0$. Using Lemma \ref{disklem}
(ii) we conclude that $\Ds$ is the minimal exhausting
disk.

Now suppose that $\lambda<p/(p-1)$. Then $\bar{w}$ satisfies the equation
\[
   \bar{w}^p = \bar{g},
\]
where $\bar{g}=\bar{g}_kt_1^k+\ldots\bar{g}_mt_1^m$ is a polynomial of degree
$m$, divisible by $t_1^k$. The affine curve $W^\circ=\Spec k[t_1,\bar{w}]$ is
irreducible and has a unibranched singularity precisely over each point of
$Z^\circ=\Spec k[t_1]$ where the differential $d\bar{g}$ has a zero. Since
$1<k<m$ and $(p,k)=(p,m)=1$ it follows that $W^\circ$ has at least two
singular points. Using Lemma \ref{disklem} (ii) again we conclude that $\Ds$
is the minimal exhausting disk.
\Endproof



\subsection{}  \label{solv5}

We can now prove the general case of Proposition \ref{solvprop}. By the result
of the previous subsections we may assume that $G$ has a proper normal
subgroup $H\lhd G$. Set $\Zs:=\Ys/H\to\Xs$ and consider the factorization
\[
   \Ys\lpfeil{H}\Zs\lpfeil{G/H}\Xs
\]
of $\phi$ into regular Galois subcovers. 

Suppose $\Zs$ is not a disk. Then by the induction hypothesis there exists a
maximal boundary component $\Xs_1\subset\Xs$ which is separating with respect
to $\Zs\to\Xs$.  If $\Zs$ is an open disk then we set $\Xs_1:=\Xs$.

Choose a connected component $\Ys_1$ of the inverse image of $\Xs_1$ in
$\Ys$. Let $G_1\subset G$ denote the stabilizer of $\Ys_1$ in $G$ and set
$H_1:=G_1\cap H_1$. Note that $H_1\lhd G_1$ is a normal subgroup.  We can
identify the quotient $\Zs_1:=\Ys_1/H_1$ (resp.\ $\Ys_1/G_1$) with a
connected component of the inverse image of $\Xs_1$ in $\Zs$ (resp.\ with
$\Xs_1$). By construction, $\Zs_1$ is either an open disk (and then
$\Zs_1=\Zs$) or an open annulus. Applying the induction hypothesis once more
to the $H_1$-cover $\Ys_1\to\Zs_1$ we obtain a maximal separating boundary
domain $\Zs_2\subset\Zs_1$ with respect to $\Ys_1\to\Zs_1$.

We claim that the subset $\Zs_2\subset\Zs_1$ is fixed by the action of
$G_1/H_1$. Indeed, any element $g\in G_1$ induces an isomorphism of the cover
$\Ys_1\to\Zs_1$. It follows that $g(\Zs_2)\subset\Zs_1$ is also a maximal
separating boundary domain with respect to $\Ys_1\to\Zs_1$. Uniqueness shows
that $g(\Zs_2)=\Zs_2$. 

The quotient $\Xs_2:=\Zs_2/(G_1/H_1)$ can be identified with a
boundary domain of $\Xs$. We claim that $\Xs_2$ is separating with respect to
the cover $\phi:\Ys\to\Xs$, and is maximal with respect to this property. 
Indeed, let $\Ys_2$ denote a connected component of the inverse image of
$\Xs_2$ in $\Ys$. We may assume that $\Ys_2$ is contained in $\Ys_1$. Then
$\Ys_2$ is also a connected component of the inverse image of
$\Zs_2\subset\Zs_1$ in $\Ys_1$. By the choice of $\Zs_2$ this means that
$\Ys_2$ is an open annulus. This shows that $\Xs_2$ is separating with respect
to the cover $\phi$. The maximality of $\Xs_2$ is proved in a similar manner. 
This completes the proof of Proposition \ref{solvprop}.
\Endproof

\section{The nonsolvable case} \label{nonsolv}

\subsection{} \label{nonsolv1}

We return to the situation considered in \S \ref{diskcover}: we are given a
regular $G$-Galois cover $\phi:\Ys\to\Xs$, where $\Xs$ is an open disk and
$\Ys$ is an open analytic curve which is {\em not} an open disk. Our goal is
to show that the set $\D$ of all affinoid disks $\Ds\subset\Xs$ which are
exhausting with respect to $\phi$ has a unique minimal element (Theorem
\ref{diskthm}). In view of Proposition \ref{solvprop} we may assume that the
group $G$ is {\em not} solvable.

\subsection{} \label{nonsolv2}

We let $\X=\Spf(A)$ and $\Y=\Spf(B)$ denote the canonical formal
models.
We let $\eta\in\partial\Xs$ denote the unique boundary point
of $\Xs$. We let $k[\eta]\subset k(\eta)$ denote the valuation ring of the
residue field of $\eta$. In this section, we consider $\eta$ as a morphism of
formal schemes $\eta:\Spf(k[\eta])\to\X$.

Let $\Ds\subset\Xs$ be an affinoid disk. By a local variant of the procedure
described in \S \ref{exhaust}, $\Ds$ gives rise to a diagram of formal
$R$-schemes
\[\begin{CD}
      \Y' @>>> \Y \\
      @VVV     @VVV \\
      \X' @>>> \X, \\
\end{CD}\]
as follows. Let $t\in A$ be a parameter and $\epsilon\in\abs{K^\times}$ such
that $\Ds=\Xs(\abs{t}\leq\epsilon)$. Let $\X'\to\X$ be the formal blowup of
the ideal $I:=(t,a)\lhd A$. Let $Z\subset\X'$ be the exceptional fiber (it is
equal to the reduced subscheme $(\X)^{\rm red}$, and it is isomorphic to
$\PP^1_k$). The morphism $\xi:\Spf(k[\eta])\to\X$ lifts uniquely to a
morphism $\xi':\Spf(k[\eta])\to\X'$.  Let $z\in Z$ denote the image of
$\xi'$ and $Z^\circ:=Z\backslash\{z\}$. Then
\[
      \Ds =]Z^\circ[_{\X'}.
\]

Let $\Y'$ be the normalization of the formal scheme $\X'$ in $\Ys$ (see
??). We call $\Y'$ the formal model of $\Ys$ induced by $\Ds$. Let
$W:=(\Y')^{\rm red}$ denote the reduced subscheme. Note that $W$ is a
connected projective $k$-curve. The canonical morphism $\Y'\to\X'$ restricts
to a finite map $W\to Z$. Let $\partial W\subset W$ denote the inverse image
of $z$ and $W^\circ:=W\backslash\partial W$. We have
\[
     \phi^{-1}(\Ds) = ]W^\circ[_{\Y'}.
\]
It follows that $\Ds$ is exhausting with respect to $\phi$ if and only if the
residue classes $]w[_{\Y'}$ are open annuli, for all $w\in\partial
W$. Actually, since $G$ acts transitively on the set $\partial W$, it suffices
that  this holds for one $w\in\partial W$.

Given a boundary point $\xi\in\partial\Ys$ the morphism
$\xi:\Spf(k[\xi])\to\Y$ lifts uniquely to $\xi':\Spf(k[\xi])\to\Y'$, and
the image $\bar{\xi}\in W$ of $\xi'$ lies in $\partial W$. We obtain a
surjective $G$-equivariant map
\begin{equation} \label{formalboundarymapeq}
    \partial\Ys\to\partial W,
\end{equation}
the local analog of the map \eqref{boundarymapeq}. If $\Ds$ is exhausting,
then this map is a bijection.

\subsection{} \label{nonsolv3}

We fix a boundary point $\xi\in\partial\Ys$. Let $\Ds\subset\Xs$ be
either an affinoid disk or the empty set. Then we let $\Vs(\Ds)$ denote the
connected component of $\phi^{-1}(\Xs\backslash\Ds)$ which `contains'
$\xi$. This means the following. If $\Ds=\emptyset$ then $\Vs(\Ds)$ is equal
to $\Ys$. On the other hand, if $\Ds$ is an affinoid disk then
$\Vs(\Ds):=]w[_{Y_R'}$, where $w\in\partial W$ is the image of $\xi$ under the
map \eqref{formalboundarymapeq}. Note that $\Ds\in\D$ if and and only if
$\Vs(\Ds)$ is an open annulus.

We let $G(\Ds)\subset G$ denote the stabilizer in $G$ of the component
$\Vs(\Ds)$. Note that $\Vs(\Ds)/G(\Ds)\cong \Xs\backslash\Ds$.

\begin{lem} \label{nonsolvlem1}
  For $\Ds\in\D$ the following holds.
  \begin{enumerate}
  \item The group $G(\Ds)$ is equal to the stabilizer $G_\xi$ in $G$ of the
    boundary point $\xi$ and is solvable. 
  \item Let $\Ds'\subset \Ds$ be a subset which is either empty or an affinoid
    disk strictly contained in $\Ds$. Then $\Vs(\Ds')\backslash\phi^{-1}(\Ds)$
    is absolutely connected. In particular, $\phi^{-1}(\Ds)$ is absolutely
    connected.
  \end{enumerate}
\end{lem}

\proof Since $\Ds\in\D$, the map \eqref{formalboundarymapeq} is a
bijection. The equality $G(\Ds)=G_\xi$ follows immediately. Moreover,
$G_\xi\subset G=\Gal(\Frac(B)/\Frac(A))$ is the decomposition group of the
discrete valuation on $\Frac(B)$ corresponding to $\xi$. We obtain a short
exact sequence
\[
      1 \to I_\xi \to G_\xi \to \Gal(k(\xi)/k(\eta)) \to 1,
\]
where $I_\xi$ is the inertia group of $\xi$. The residue field extension
$k(\xi)/k(\eta)$ is a finite extension of complete discrete valued fields with
algebraic residue field. So $\Gal(k(\xi)/k(\eta))$ is also an inertia
group and hence solvable. We conclude that $G_\xi$ is solvable.

For the proof of (ii) we first assume that $\Ds'=\emptyset$, and we use the
notation introduced in \S \ref{nonsolv2}. We have already noted that $W$ is
connected. Since $\Ds\in\D$, the subset $\partial W$ consists of
smooth points of $W$. It follows that the complement
$W^\circ=W\backslash\partial W$ is still connected. Now Remark \ref{rigidrem}
(iv),(v) shows that
\[
   \phi^{-1}(\Ds) = ]W^\circ[_{\Y'}
\]
is absolutely connected, proving (ii) if $\Ds=\emptyset$. The proof in the
case $\Ds'\neq \emptyset$ is similar and left to the reader.
\Endproof

\subsection{}

Let $\Xs\an$ denote the Berkovich analytic space associated to $\Xs$, see
\cite{Berkovich90}. As a set, $\Xs\an$ consists of all continuous multiplicative
seminorms $\abs{\,\cdot\,}_x:A\to\RR_{\geq 0}$ bounded by $1$ which extend the
standard valuation $\abs{\,\cdot\,}$ on $R$. To each point
$\abs{\,\cdot\,}_x\in\Xs\an$ we can associate its residue field $\HH(x)$,
which is defined as the completion of the fraction field of $A/{\rm
  Ker}(\abs{\,\cdot\,}_x)$. By construction, $\HH(x)/K$ is an extension of
complete valued fields. We let $\widetilde{\HH(x)}$ denote the residue field
of $\HH(x)$.

Any point $x\in\Xs$ gives rise to a point $\Xs\an$ by the formula
$\abs{f}_x:=\abs{f(x)}$. We may thus consider $\Xs$ as a subset of $\Xs\an$
(called the set of {\em classical points}). Classical points are characterized
by the property that the extension $\HH(x)/K$ is finite.

In order to have a uniform and suggestive notation, we shall write
$x\in\Xs\an$ instead of $\abs{\,\cdot\,}_x\in\Xs\an$ and $\abs{f(x)}$ instead
of $\abs{f}_x$, for arbitrary points on $\Xs\an$. For instance, to any closed
disk $\Ds\subset\Xs$ (affinoid or not) we can associate a point
$x_\Ds\in\Xs\an$ by setting
\[
     \abs{f(x_\Ds)} := \max_{x'\in\Ds} \abs{f(x')}.
\]
If $\Ds$ is affinoid, then the residue field $\widetilde{\HH(x)}$ can be
identified with the function field of the canonical reduction of $\Ds$. In
particular, $\widetilde{\HH(x)}/k$ has transcendence degree one. Otherwise,
$\widetilde{\HH(x)}=k$. See \cite{Berkovich90}, \S 1.4.4.

The space $\Xs\an$ carries a natural topology which makes it a locally compact
Hausdorff space. If $\Ds\in\Xs$ is an affinoid disk, then
$\Ds\an\subset\Xs\an$ is a compact subset. It follows that the limit
\[
      x:=\lim_{\Ds\in\D} x_\Ds \in\Xs\an
\]
exists. In fact, it is easy to see that for all $f\in A$ we have
\[
      \abs{f(x)} = \inf_{\Ds\in\D} \abs{f(x_\Ds)}.
\]

\begin{prop} \label{limdiskprop}
  Let $\Ds_0:=\cap_{\Ds\in\D}\Ds$. Exactly one of the following cases occurs.
  \begin{itemize}
  \item[(1)]
    The limit point $x$ is a classical point, and $\Ds_0=\{x\}$. 
  \item[(2)] The set $\Ds_0$ is an affinoid disk, and $x=x_{\Ds_0}$. 
  \item[(3)] The set $\Ds_0$ is a closed disk which is not affinoid, and
    $x=x_{\Ds_0}$. 
  \item[(4)]
    The set $\Ds_0$ is empty. 
  \end{itemize}
  In Case (1), (3) and (4) we have $\widetilde{\HH(x)}=k$.
\end{prop}

\proof This follows from the classification of points on $(\AA_K^1)\an$ in
\cite{Berkovich90}, \S 1.4.4.
\Endproof

\begin{lem} \label{nonsolvlem2}
  Assume that we are in Case (1), (3) or (4) of Proposition
  \ref{limdiskprop}. Then for any $y\in\phi^{-1}(x)$, the stabilizer
  $G_y\subset G$ of $y$ is solvable.
\end{lem}

\proof Let $\OO_{\Xs\an,x}$ denote the local ring of the point $x$ on the
analytic $K$-space $\Xs\an$. By \cite{Berkovich93}, Theorem 2.1.5,
$\OO_{\Xs\an,x}$ is a henselian local ring. Moreover, by {\em loc.cit.}, Lemma
2.1.6, we have a decomposition
\[
     \phi_*(\OO_{\Ys\an})_x = \prod_{i=1}^n \OO_{\Ys\an,y_i},
\]
where $\phi^{-1}(x)=\{y_1,\ldots,y_n\}$. Therefore, the extension
$\OO_{\Ys\an,y}/\OO_{\Xs\an,x}$ is a Galois extension of henselian local rings
with Galois group $G_y$. It follows that $G_y$ sits in a short exact sequence
\[
      1 \to I_y \to G_y \to \Gal(\kappa(y)/\kappa(x)) \to 1,
\]
where $\kappa(x)$ and $\kappa(y)$ are the residue fields of $\OO_{\Xs\an,x}$
and $\OO_{\Ys\an,y}$, respectively, and $I_y$ is the inertia group.  By
\cite{Berkovich93}, Proposition 2.4.3 and Proposition 2.4.4, $I_y$ is
solvable. It remains to show that $\Gal(\kappa(y)/\kappa(x)$ is solvable. By
\cite{Berkovich93}, Theorem 2.3.3, $\kappa(x)$ is a henselian valued field
(called {\em quasi-complete} in {\em loc.cit.}) whose completion is the field
$\HH(x)$. By Proposition \ref{limdiskprop}, the residue field of $\HH(x)$
(equal to the residue field of $\kappa(x)$) is equal to $k$, which is
algebraically closed by assumption. Using again \cite{Berkovich93},
Proposition 2.4.4 we conclude that $\Gal(\kappa(y)/\kappa(x))$ and hence $G_y$
is solvable.
\Endproof

\begin{lem} \label{nonsolvlem3}
\begin{enumerate}
\item In Case (1) and (4) of Proposition \ref{limdiskprop}, the affinoid
  disks $\Ds\an$, $\Ds\in\D$, form a neighborhood basis of $x$ in $\Xs\an$.
\item
  In Case (3), a neighborhood basis of $x$ in $\Xs\an$ is given by the sets
  $(\Ds\backslash\Ds')\an$, where $\Ds\in\D$ and $\Ds'$ is an affinoid disk
  contained in $\Ds_0=\cap_{\Ds\in\D}\Ds$. 
\end{enumerate}
\end{lem}

\proof In Case (1), the statement in (i) is clear. Suppose we are in Case (4),
i.e.\ $\cap_{\Ds\in\D}\Ds=\emptyset$. Any $f\in A$, seen as an analytic
function on $\Xs\an$, has finitely many zeroes, all of which are classical
points. Therefore, there exists $\Ds_1\in\D$ such that $f|_{\Ds\an}$ is an
invertible analytic function on $\Ds\an$, for all $\Ds\in\D$ with
$\Ds\subset\Ds_1$. Applying the maximum principle to $f|_{\Ds\an}$ and
$f^{-1}|_{\Ds\an}$ we see that $\abs{f(x)}=\abs{f(x_\Ds)}$. The statement in
(i) now follows from the definition of the topology of $\Xs\an$. The proof of
(ii) is similar.
\Endproof

\subsection{} \label{nonsolv4}

We can now finish the proof of Theorem \ref{diskthm}.  We first suppose that
we are in Case (1) or (4) of Proposition \ref{limdiskprop}. Let
$\phi^{-1}(x)=\{y_1,\ldots,y_n\}$ be the fiber above $x$. Then it follows from
Lemma \ref{nonsolvlem3} and \cite{Berkovich93}, proof of Theorem 2.1.5, that
for all sufficiently small $\Ds\in\D$ the inverse image $\phi^{-1}(\Ds)$
decomposes into $n$ disjoint affinoid neighborhoods of the points $y_i$. But
$\phi^{-1}(\Ds)$ is connected by Lemma \ref{nonsolvlem1} (ii), hence $n=1$. By
Lemma \ref{nonsolvlem2} this shows that the group $G$ is solvable,
contradicting our assumption. We conclude that Case (1) and (4) of Proposition
\ref{limdiskprop} cannot occur.

We now suppose that we are in Case (2) of Proposition \ref{limdiskprop}, i.e.\
$\Ds_0:=\cap_{\Ds\in\D}\Ds$ is an affinoid disk. By Lemma \ref{nonsolvlem1}
(i) we have
\[
     g(\Vs(\Ds))\cap\Vs(\Ds) = \emptyset
\]
for all $\Ds\in\D$ and $g\in G\backslash G_\xi$. It follows easily that 
\[
      \Vs(\Ds_0)=\cup_{\Ds\in\D} \Vs(\Ds) \quad\text{and}\quad
        G(\Ds_0)=G_\xi.
\]
Consider the $G_\xi$-cover
\[
      \Vs(\Ds_0)\to \Xs\backslash\Ds_0.
\]
If $\Ds\subset\Ds_0$ is an affinoid disk, then $\Ds\in\D$ if and only if
$\Us:=\Xs\backslash\Ds$ is a separating boundary component for this cover, see
\S \ref{solv}.  Since $G_\xi$ is solvable by Lemma \ref{nonsolvlem1} (i),
Proposition \ref{solvprop} shows that there exists a unique maximal
separating boundary component. It follows that $\Ds_0\in\D$ is the unique
maximal element of $\Ds$, i.e.\ Theorem \ref{diskthm} holds.

It remains to rule out Case (3) of Proposition \ref{limdiskprop}. Using Lemma
\ref{nonsolvlem3} (ii) we find an element $\Ds\in\D$ and an affinoid disk
$\Ds'\subset\Ds$ such that $x\in\Ds\an\backslash(\Ds')\an$ and such that
$\Vs(\Ds')\an\backslash\Vs(\Ds)\an$ contains a unique point
$y\in\phi^{-1}(x)$. By Lemma \ref{nonsolvlem2} this implies that the
stabilizer of $\Vs(\Ds')\backslash\Vs(\Ds)$ in $G$ is solvable. Together with
Lemma \ref{nonsolvlem1} (ii) this shows that $G(\Ds')$ is solvable. Applying
Proposition \ref{solvprop} to the $G(\Ds')$-cover
\[
      \Vs(\Ds')\to\Xs\backslash\Ds'
\]
and using a similar argument as in Case (2) we conclude that $\D$ has a
minimal element. But since this minimal element must be equal to the
intersection $\Ds_0=\cap_{\Ds\in\D}\Ds$ which is {\em not} an affinoid disk in
Case (3), we obtain a contradiction. Therefore, Case (3) of Proposition
\ref{limdiskprop} cannot occur. This completes the proof of Theorem
\ref{diskthm}.
\Endproof

\bibliographystyle{plain} \bibliography{references}

\begin{thebibliography}{10}

\bibitem{Abbes2000}
A.~Abbes.
\newblock {R\'eduction semi-stable des courbes d'apres Artin, Deligne,
  Grothendieck, Mumford, Saito, Winters,..}
\newblock In {\em Courbes semi-stables et groupes fondamental en g\'eometrie
  alg\'ebrique}, pages 59--110. Birkh\"auser, 2000.

\bibitem{KaiDiss}
K.~Arzdorf.
\newblock {\em Semistable reduction of cyclic covers of prime power degree}.
\newblock PhD thesis, Leibniz University Hannover, 2012.

\bibitem{Berkovich90}
V.G. Berkovich.
\newblock {\em Spectral theory and analytic geometry over non-archimedian
  fields}.
\newblock Number~33 in Mathematical Surveys and Monographs. AMS, 1990.

\bibitem{Berkovich93}
V.G. Berkovich.
\newblock {\'Etale cohomology of non-archimedian analytic spaces}.
\newblock {\em Publ.\ Math.\ IHES}, 78:5--161, 1993.

\bibitem{Berthelot96}
P.~Berthelot.
\newblock Cohomologie rigide et cohomologie rigid \`a supports propres.
\newblock Prepublication 96-03, Universit\'e de Rennes 1.

\bibitem{BGR}
S.~Bosch, U.~G\"untzer, and R.~Remmert.
\newblock {\em {Non-Archimedian analysis: a systematic approach to rigid
  analytic geometry}}.
\newblock Number 261 in Grundlehren der math.\ Wiss. Springer-Verlag, 1984.

\bibitem{BoschLuetkebohmert85}
S.~Bosch and W.~L\"utkebohmert.
\newblock {Stable reduction and uniformization of abelian varieties, I}.
\newblock {\em Math.\ Ann.}, 270(3):349--379, 1985.

\bibitem{BoschLuetkebohmert93}
S.~Bosch and W.~L\"utkebohmert.
\newblock {Formal and rigid gemetry: I. Rigid spaces}.
\newblock {\em Math.\ Ann.}, 295:291--317, 1993.

\bibitem{Conrad99}
B.~Conrad.
\newblock Irreducible components of rigid spaces.
\newblock {\em Annales Inst.\ Fourier}, 49(2):473--541, 1999.

\bibitem{deJong95}
A.J. deJong.
\newblock {Crystalline Dieudonn\'e module theory via formal and rigid
  geometry}.
\newblock {\em Publ.\ Math.\ IHES}, 82:5--96, 1995.

\bibitem{DeligneMumford69}
P.\ Deligne and D.\ Mumford.
\newblock The irreducibility of the space of curves of given genus.
\newblock {\em Publ.\ Math.\ IHES}, 36:75--109, 1969.

\bibitem{Epp73}
H.P.\ Epp.
\newblock Eliminating wild ramification.
\newblock {\em Inventiones math.}, 19:235--249, 1973.

\bibitem{FresnelvdPut}
J.~Fresnel and M.~van~der Put.
\newblock {\em Rigid Analytic Geometry and its Applications}.
\newblock Number 218 in Progress in Math. Birkh\"auser, 2004.

\bibitem{LehrMatignon06}
C.~Lehr and M.~Matignon.
\newblock Wild monodromy and automorphisms of curves.
\newblock {\em Duke Math.\ J.}, 135(3):569--586, 2006.

\bibitem{LiuAG}
Q.~Liu.
\newblock {\em Algebraic geometry and arithmetic curves}.
\newblock Oxford Univ.\ Press, 2002.

\bibitem{Matignon03}
M.~Matignon.
\newblock Vers un algorithme pour la r\'eduction semistable des rev\^etements
  $p$-cycliques de la droite projective sur un corps $p$-adique.
\newblock {\em Math.\ Ann.}, 325(2):323--354, 2003.

\bibitem{RaynaudLecture1}
M.~Raynaud.
\newblock {Enonc\'es de permanence en g\'eom\'etrie relative}.
\newblock {Journ\'ees de G\'eom\'etrie Arithm\'etiques de Rennes}, July 2009.

\bibitem{Temkin10}
M.~Temkin.
\newblock Stable modification of relative curves.
\newblock {\em J.\ Algebraic Geometry}, 19:603--677, 2010.

\end{thebibliography}

\end{document}